\documentclass[11pt]{article}

\usepackage{epsfig,graphics}
\usepackage{amssymb,amsmath,amsthm,bm}

\newtheorem{theo}{Theorem}

\newtheorem{cor}{Corollary}
\newtheorem{mrem}{Remark}
\newtheorem{mdef}{Definition}

\def\R{{\mathbb R}}

\def\RR{{\mathbb R}}

\def\ga{\alpha}

\def\gth{\theta}

\def\gb{\beta}

\def\gs{\sigma}
\def\gl{\lambda}
\def\wt{\widetilde}

\def\gp{{\prime}}

\def\ep{\epsilon}
\def\vep{\varepsilon}

\def\gt{\triangle}

\def\b0{{\bf 0}}
\def\1{{\bf 1}}

\def\bd{{\bf d}}

\def\cB{\mathcal B}
\def\cD{\mathcal D}

\def\cG{\mathcal G}
\def\cH{\mathcal H}

\def\rem{{\rm rem}}

\def\err{{\rm res}}

\def\terr{{\rm err}}
\def\tErr{{\rm Err}}

\def\bd{{\rm bd}}
\def\Bd{{\rm Bd}}

\def\wh{\widehat}

\def\wt{\widetilde}

\makeatletter
\def\widebreve{\mathpalette\wide@breve}
\def\wide@breve#1#2{\sbox\z@{$#1#2$}%
     \mathop{\vbox{\m@th\ialign{##\crcr
\kern0.08em\brevefill#1{0.8\wd\z@}\crcr\noalign{\nointerlineskip}%
                    $\hss#1#2\hss$\crcr}}}\limits}
\def\brevefill#1#2{$\m@th\sbox\tw@{$#1($}%
  \hss\resizebox{#2}{\wd\tw@}{\rotatebox[origin=c]{90}{\upshape(}}\hss$}
\makeatletter
\def\wb{\widebreve}

\usepackage{scalerel,stackengine}
\stackMath
\newcommand\reallywidecheck[1]{%
\savestack{\tmpbox}{\stretchto{%
  \scaleto{%
    \scalerel*[\widthof{\ensuremath{#1}}]{\kern-.6pt\bigwedge\kern-.6pt}%
    {\rule[-\textheight/2]{1ex}{\textheight}}
  }{\textheight}%
}{0.5ex}}%
\stackon[1pt]{#1}{\scalebox{-1}{\tmpbox}}%
}
\def\wc{\reallywidecheck}

\begin{document}

\title{Analysis of  an  Adaptive Short-Time Fourier Transform-Based  
Multicomponent Signal Separation Method Derived from \\
Linear Chirp Local Approximation\thanks{
This work is partially supported by the Hong Kong Research Council, under Projects $\sharp$ 12300917 and $\sharp$ 12303218, and HKBU Grants $\sharp$ RC-ICRS/16-17/03 and $\sharp$ RC-FNRA-IG/18-19/SCI/01, the Simons Foundation, under Grant $\sharp$ 353185, and the National Natural Science Foundation of China 
under Grants $\sharp$ 61972265 and $\sharp$ 11871348, the Natural Science Foundation of Guangdong Province of China  under Grant $\sharp$ 2020B1515310008.}
}

\author{Charles K. Chui${}^{1}$,  Qingtang Jiang${}^2$, Lin Li${}^3$,
and  Jian Lu${}^{4}$}

\date{} 

\maketitle


\bigskip
{\small 1. 
College of Mathematics \& Statistics, Shenzhen University, Shenzhen 518060, China}

\quad {\small and Department of Mathematics, Hong Kong Baptist University, Hong Kong.}

{\small 2. Department of Math \& Computer Sci., 
Univ. of Missouri-St. Louis, St. Louis,  MO 63121, USA.} 

{\small 3. School of Electronic Engineering, Xidian University, Xi${}'$an 710071, China.}  

{\small 4. Shenzhen Key Laboratory of Advanced Machine Learning and Applications,}  

\quad {\small College of Mathematics \& Statistics, Shenzhen University, Shenzhen 518060, China.} 

\begin{abstract}
The synchrosqueezing transform (SST) has been developed as a powerful EMD-like tool for instantaneous frequency (IF) estimation and component separation of non-stationary multicomponent signals. Recently, a direct method of the time-frequency approach, called signal separation operation (SSO), was introduced to solving the problem of multicomponent signal separation. While both SST and SSO are mathematically rigorous on IF estimation, SSO avoids the second step of the two-step SST method in component recovery (mode retrieval). 
In addition, SSO is simple: the IF of a component 
is estimated by  a time-frequency ridge of the SSO plane; and this component is recovered by simply plugging the time-frequency ridge to the SSO operation. 
In recent paper \lq\lq{}{\it Direct signal separation via extraction of local frequencies with adaptive time-varying parameters}\rq\rq{},  after showing that  the SSO operation is related to the adaptive short-time Fourier transform (STFT),  the authors obtained a more accurate component recovery formula derived from the  linear chirp (also called linear frequency modulation signal) approximation at any local time and they also proposed a recovery scheme to extract the signal components one by one with the time-varying window updated for each component.  However the theoretical analysis of 
the recovery formula derived from linear chirp local approximation has not been studied there. In this paper, we carry out such analysis and obtain error bounds for 
 IF estimation and component recovery. These results provide a mathematical guarantee to the proposed adaptive STFT-based non-stationary multicomponent signal separation method.
\end{abstract}

\bigskip 

\bigskip 

\centerline{{\bf AMS Mathematics Subject Classification:} 42A38, 42C40,  42C15}


\section{Introduction}
Many real-world signals are represented as a superposition of Fourier-like oscillatory modes:
\begin{equation}
\label{AHM0}
x(t)=A_0(t)+\sum_{k=1}^K x_k(t),  \quad x_k(t)=A_k(t) \cos \big(2\pi \phi_k(t)\big), 
\end{equation}
where $A_k(t), \phi_k'(t)>0$ and $A_k(t)$ varies slowly. 
 \eqref{AHM0} is also called an {adaptive harmonic model (AHM) representation of $x(t)$, with $A_0(t)$ called the trend, $A_k(t)$ instantaneous amplitudes (IAs)  and $\phi'_k(t)$ the instantaneous frequencies (IFs).  To represent $x(t)$ as \eqref{AHM0} is important for extracting information hidden in the non-stationary source signal $x(t)$. 
 The empirical mode decomposition (EMD) scheme 
is a popular method to decompose non-stationary signals \cite{Huang98}. 
EMD 
decomposes a signal into its ``IMFs" and ``trend" without the concern of recovering the true IMFs and trend of the source signal.
In this paper we  consider recovering trend $A_0(t)$ and  components $x_k(t)$ of the source signal $x(t)$.

The continuous wavelet transform (CWT)-based  synchrosqueezed wavelet transform (SST), introduced in  \cite{Daub_Maes96} and further developed in \cite{Daub_Lu_Wu11}, 
provides mathematical theorems to guarantee the recovery of oscillatory modes.
Since the seminal work \cite{Daub_Lu_Wu11} was available to the research community, various  SSTs have been proposed and studied, see e.g. 
\cite{BMO18}-\cite{Chui_Lin_Wu15}, \cite{Chui_Walt15},
\cite{Jiang_Suter17}-\cite{ LCJ18},
 \cite{LJL20}-\cite{Yang15}. 


Recently, a direct method of the time-frequency approach, called signal separation operation or signal separation operator (SSO), was introduced in \cite{Chui_Mhaskar20} to solving the problem of multicomponent signal separation. While both SST and SSO are mathematically rigorous on IF estimation, SSO avoids the second step of the two-step SST method in component recovery, which depends heavily on the accuracy of the estimated IFs. In addition, SSO is simple: the IF $\phi_k\rq{}(t)$ of the $k$-th component $x_k(t)$ of a multicomponent signal $x(t)$ is estimated by a time-frequency ridge $\widehat \gth_k(t)$ of the SSO operation $T_x(t, \gth)$; and $x_k(t)$ is recovered by simply plugging  $\widehat \gth_k(t)$ to  $T_x(t, \gth)$: $x_k(t)\approx  T_x(t, \wh \gth_k(t))$. 

In their recent paper \cite{LCJ20}, the authors  show that  the SSO operation is related to the adaptive STFT defined by 
\begin{eqnarray}
\label{def_STFT_para1}
\wt V_{x}(t, \eta) 
\hskip -0.6cm &&:=\int_{-\infty}^\infty x(\tau)g_{\gs(t)}(\tau-t)e^{-i2\pi \eta (\tau-t)}d\tau\\
\label{def_STFT_para2}
&&=\int_{-\infty}^\infty x(t+\tau)\frac 1{\gs(t)}g (\frac \tau{\gs(t)})e^{-i2\pi \eta\tau}d\tau,
\end{eqnarray}
where $t$ and $\eta$ are the time variable and the frequency variable respectively, 
$\gs=\gs(t)$ is a positive function of $t$, and $g_{\gs(t)}(\tau)$ is defined by
\begin{equation}
\label{def_g_dilation}
g_{\gs(t)}(\tau):=\frac 1{\gs(t)}g(\frac \tau {\gs(t)}),
\end{equation}
with $g\in L_2(\R)$. The window width of $g_{\gs(t)}(\tau)$ is $\gs(t)$ (up to a constant), depending on the time variable $t$.

The authors of \cite{LCJ20} obtained a more accurate component recovery formula derived from linear chirp (also called linear frequency modulation signal) approximation at any local time and they also proposed a recovery scheme to extract the signal components one by one with the time-varying window updated for each component.  However the theoretical analysis of the recovery  formula derived from the linear chirp local approximation has not been studied there. In this paper, we carry out such analysis and obtain error bounds for 
 IF estimation and component recovery with this new formula. These results provide a mathematical guarantee to the proposed adaptive STFT-based
 non-stationary multicomponent signal separation method.

Recently our linear chirp local approximation-based approach has been extended in \cite{LHJC20} to 3-dimensional case with variables of time, frequency and chirp rate to recover components with crossover IFs.    
Very recently the STFT-based SSO method is extended to the continuous wavelet transform ( CWT)-based SSO in  \cite{Chui_Han20,Chui_Mhaskar20}.  Furthermore,  a more accurate component recovery formula based on the adaptive CWT was derived based on linear chirp local approximation in \cite{CJLL20}. 

In this paper we consider the complex version 
of multi-component signals $x(t)$ of \eqref{AHM0}:
\begin{equation}
  \label{AHM}
  x(t)=A_0(t)+\sum_{k=1}^K x_k(t)=A_0(t)+\sum_{k=1}^K A_k(t) e^{i2\pi\phi_k(t)}
  \end{equation}
 with $A_k(t), \phi_k'(t)>0$. For convenience, we let $\phi_0(0)=0$. Thus the trend can also be written as  $x_0(t):=A_0(t) e^{i2\pi\phi_0(t)}$.   The main goal of this paper is to analyze the error bound for the recovery formula based on linear chirp approximation (called the linear chirp-based model in this paper). To make this paper be self-contained, we also include the error bound for the recovery formula based on sinusoidal signal approximation (called 
 the sinusoidal signal-based model in this paper).  The sinusoidal signal-based model and 
 the linear chirp-based model will be studied in Sections 2 and 3 respectively. 

Note that when $\gs(t) \equiv\gs$, a positive constant, then 
$\wt V_{x}(t, \eta)$ in \eqref{def_g_dilation} is the conventional (modified) STFT of $x(t)$ with window function 
$h(t)=\frac 1\gs g(\frac t\gs)$: 
\begin{eqnarray}
\label{def_STFT}
V_x(t, \eta) \hskip -0.6cm&&:=\int_{-\infty}^\infty x(\tau) h(\tau-t)e^{-i2\pi \eta (\tau-t)}d\tau.
\end{eqnarray}
Thus our linear chirp local approximation-based recovery formula and its analysis established in Section 3 apply to  the conventional STFT-based approach. 

\section{Sinusoidal signal-based model}

Let $x(t)$ be a non-stationary multicomponent signal of the form \eqref{AHM}. 
We assume $A_k(t), \phi_k(t)$ satisfy  
\begin{eqnarray}
 \label{cond_basic0}&&A_k(t)\in C^1(\R)\cap L_\infty(\R), \phi_k(t)\in C^2(\R), \\
\label{cond_basic}&& A_k(t)>0, \; \inf_{t\in \R} \phi^\gp_k(t)>0, \; \sup_{t\in \R} \phi^\gp_k(t)<\infty, \\
\label{freq_resolution_adp}
&& d\rq{}:=\min_{k\in \{2, \cdots, K\}}\min_{t\in \R}(\phi^\gp_k(t)-\phi^\gp_{k-1}(t))> 0.
\end{eqnarray}
In this section we consider the case that each component $x_k(t)=A_k(t)e^{i2\pi \phi_k(t)}$ is approximated locally by a sinusoidal signal. More precisely, we assume $A_k(t)$ 
and $\phi^\gp_k(t)$ changes slowly:   
\begin{eqnarray}
\label{condition1} && |A_k(t+\tau)-A_k(t)|\le \vep_1 |\tau|\;  A_k(t),  \; t\in \R, \; 0\le k\le K, 
\\ 
\label{condition1b}&& |\phi_k^{\gp\gp}(t)|\le \vep_2, \; t\in \R, \; 1\le k\le K, 
\end{eqnarray}
for some positive number $\vep_1, \vep_2$. Let ${\cB}_{\vep_1, \vep_2}$ denote the set of multicomponent signals of \eqref{AHM} satisfying \eqref{cond_basic0}-\eqref{condition1b}. 

Denote 
\begin{equation}\label{def_M_u}
\mu=\mu(t):= \min_{0\leq k\leq K} A_{k}(t), ~~ M=M(t):=\sum_{k=0}^K A_{k}(t), ~~M_\ell=M_\ell(t):=\sum_{k\not= \ell} A_{k}(t).  
\end{equation}
Throughout this paper, $\sum_{k\not= \ell}$ denotes $\sum_{\{k: ~ k\not= \ell, 0\le k\le K\}}$.    

	\bigskip 

Let $x(t)\in {\cB}_{\vep_1, \vep_2}$.   Write $x_k(t+\tau)$ as 
\begin{equation}
\label{approx_sinusoidal} 
x_k(t+\tau)=x_k(t)e^{i2\pi \phi^\gp_k(t)\tau }+ x_{k, {\rm rs}}(t, \tau), 
\end{equation}
where $x_{k, {\rm rs}}(t, \tau):=x_k(t+\tau)-x_k(t)e^{i2\pi \phi^\gp_k(t)\tau }$. Note that 
for any given $t$, $x_k(t)e^{i2\pi \phi^\gp_k(t)\tau }$ (as a function of $\tau$) is a   sinusoidal function in $\tau$. We use $x_k(t)e^{i2\pi \phi^\gp_k(t)\tau }$ to approximate $x(t+\tau)$ at a local time $t$. 
By \eqref{approx_sinusoidal}, we have 
\begin{eqnarray}
\nonumber \wt V_{x}(t, \eta)
\hskip -0.6cm &&=\int_\R \sum_{k=0}^K \Big(x_k(t)e^{i2\pi \phi^\gp_k(t)\tau }+ x_{k, {\rm rs}}(t, \tau)\Big)
\frac 1{\gs(t)}g (\frac \tau{\gs(t)})e^{-i2\pi \eta\tau}d\tau\\
\label{STFT_approx_1st}
&&=\sum_{k=0}^K x_k(t)  \wh g\big(\gs(t)(\eta-\phi_k^\gp(t) \big) +\rem_0,  
\end{eqnarray}
where $\rem_0$ is the remainder for the expansion of  $\wt V_{x}(t, \eta)$  
in \eqref{STFT_approx_1st} given by 
\begin{equation*}
\rem_0:=\sum_{k=0}^K \int_\R x_{k, {\rm rs}}(t, \tau) \frac 1{\gs(t)}g(\frac \tau{\gs(t)}) e^{-i2\pi \eta \tau}d\tau. 
\end{equation*}
Writing $x_{k, {\rm rs}}(t, \tau)$ as 
$$
x_{k, {\rm rs}}(t, \tau)=(A_k(t+\tau)-A_k(t))e^{i2\pi \phi_k(t+\tau)} +x_k(t)e^{i2\pi \phi^\gp_k(t)\tau}\big(
e^{i2\pi (\phi_k(t+\tau)-\phi_k(t)-\phi_k^\gp(t) \tau)}-1\big).
$$
From \eqref{condition1b}, we have 
$$
|e^{i2\pi (\phi_k(t+\tau)-\phi_k(t)-\phi_k^\gp(t) \tau)}-1|\le 2\pi |\phi_k(t+\tau)-\phi_k(t)-\phi_k^\gp(t) \tau|
\le 
\pi\vep_2 |\tau|^2. 
$$
This, together with \eqref{condition1}, leads to 
\begin{eqnarray}
\nonumber |\rem_0|\hskip -0.6cm && \le \sum_{k=0}^K \int_\R  \vep_1 | \tau| A_k(t) \frac 1{\gs(t)}|g(\frac \tau{\gs(t)})| d\tau+ \sum_{k=1}^K A_k(t)\int_\R  \pi \vep_2 | \tau|^2 \frac 1{\gs(t)}|g(\frac \tau{\gs(t)})|
 d\tau   \\
 \label{rem0_est}
 &&=M(t)\big(\vep_1 I_1 \gs(t)  +\pi  \vep_2 I_2 \gs^2(t)\big)=:M(t)\lambda_0(t),  
\end{eqnarray}
where $I_n$ and $\gl_0(t)$ are defined by 
\begin{eqnarray}
&&\label{def_In}
\label{def_tIn}
I_n:=\int_\R  | \tau^n g(\tau)| d\tau, \\ 
&&\label{def_Lam0}
\lambda_0(t):=\vep_1  I_1\gs (t) +\pi  \vep_2 I_2 \gs^2 (t). 
\end{eqnarray}

\bigskip 

If the remainder $\rem_0$ in \eqref{STFT_approx_1st} is small, then 
the term $x_k(t)\wh g\big(\gs(t)(\eta-\phi_k^\gp(t) \big)$  in \eqref{STFT_approx_1st} governs the time-frequency zone of the STFT $\wt V_{x_k}$ of the $k$th component $x_k(t)$ of $x(t)$. In particular, if $g$ is band-limited, that is $\wh g$ is compactly supported, to say supp($\wh g)\subset [-\gt, \gt]$ for some $\gt>0$, then $x_k(t)\wh g\big(\gs(t)(\eta-\phi_k^\gp(t) \big)$ lies  within the time-frequency zone:
 $$
\{(t, \eta): |\eta-\phi_k^\gp(t)|< \frac {\gt}{\gs(t)}, t\in \R\}. 
 $$
 
 If $\wh g$ is not compactly supported,  
we need to define the {\it essential support} of  $\wh g$ outside which $\wh g(\xi)\approx 0$. 
More precisely, for a given threshold $0<\tau_0<1$, if
$|\wh g(\xi)|\le \tau_0$ for $|\xi|\ge \ga$, then we say $\wh g(\xi)$ is {\it essentially supported} in $[-\ga, \ga]$. 
When $|\wh g(\xi)|$ is even and decreasing for $\xi\ge 0$,   
then $\ga$ can be obtained by solving 
\begin{equation}
\label{def_ga_general}
|\wh g(\ga)|=\tau_0. 
\end{equation}
For example, when  $g$ is the Gaussian window function given by 
\begin{equation}
\label{def_g}
g(t)=\frac 1{\sqrt {2\pi}} e^{-\frac {t^2}2},
\end{equation}
then, with $\wh g(\xi)=e^{-2\pi^2 \xi^2}$,
the corresponding $\ga$ is given by  
\begin{equation}
\label{def_ga}
\ga=\frac 1{2\pi}\sqrt{2\ln (\frac 1{\tau_0 })} . 
\end{equation}
In the following we assume $|\wh g(\xi)|$ is even. In addition, we assume $|\wh g(\xi)|$
decreasing for $\xi\ge 0$ if $\wh g(\xi)$ is not compactly supported, and we let $|\wh g|^{-1}$ denote the inverse function of $|\wh g(\xi)|, \xi\in (0, \infty)$. If supp$ \wh g=[-\gb, \gb]$ for some $\gb>0$, then we assume 
 $|\wh g(\xi)|$ is decreasing for $\xi\in [0, \gb]$, and let $|\wh g|^{-1}$ denote the inverse function of $|\wh g(\xi)|, \xi\in (0, \gb)$. Moreover,  we assume 
$$
\int_\R g(t) dt =1. 
$$

For $g$ with $\wh g(\xi)$  essentially supported in $[-\ga, \ga]$,  we then define the time-frequency zone $Z_k$ of the $k$th-component $x_k(t)$ of $x(t)$ by 
\begin{equation}
\label{def_Zk}
Z_k:=\{(t, \eta): |\wh g\big(\gs(t)(\eta-\phi_k^\gp(t) \big) |>\tau_0, t\in \R\}=
\{(t, \eta): |\eta-\phi_k^\gp(t)|< \frac {\ga}{\gs(t)}, t\in \R\}. 
\end{equation}
Thus the multicomponent signal $x(t)$ is well-separated, if $\gs(t)$ satisfies 
\begin{equation}
\label{separated_cond_1st}
\gs(t)\ge \frac {2\ga}{ \phi_k^\gp(t)-\phi_{k-1}^\gp(t)}, \; t\in \R, k=1, 2, \cdots, K. 
\end{equation}
In this case $Z_k\cap Z_{\ell}=\emptyset, k\not=\ell$. 
In this section we assume that \eqref{separated_cond_1st} holds for some $\gs(t)$. Due to \eqref{freq_resolution_adp}, there always exists 
a bounded $\gs(t)$ such that \eqref{separated_cond_1st} holds. 
When $g$ is the Gaussian window function, for the sinusoidal signal-based adaptive method,  we may choose $\gs(t)$ 
as in \cite{LCHJJ18} to be 
 $$
 \gs_1(t):=\max_{2\le k\le K}\big\{\frac {2\ga}{ \phi_k^\gp(t)-\phi_{k-1}^\gp(t)}\big\}=\frac {2\ga}{\min_{2\le k\le K}  \{\phi_k^\gp(t)-\phi_{k-1}^\gp(t)\}}.
 $$
The reason is that from the error bound for the recovery formula (see Corollary \ref{cor:sinu_cor2} below), the smaller $\gs(t)$ is, the smaller the error bound.   

From \eqref{separated_cond_1st}, we have 
\begin{equation}
\label{different_IF_est}
\gs(t)| \phi_k^\gp(t)-\phi_\ell^\gp(t)|\ge{2\ga |k-\ell|}.  
\end{equation}
and hence for $(t, \eta)\in Z_k$, 
\begin{equation}
\label{zone_est}
|\eta- \phi_\ell^\gp(t)|\ge | \phi_\ell^\gp(t)-\phi_{k}^\gp(t)|-|\eta-\phi_{k}^\gp(t)|\ge 
 \frac {\ga (2|\ell-k|-1)}{\gs(t)}. 
\end{equation}
 
 Denote 
 \begin{equation}
 \label{def_err}
 \terr_\ell=\terr_\ell(t):=M(t)\gl_0(t)+\sum_{k\ne \ell}A_k(t) |\wh g\big(\ga (2|\ell-k|-1\big)|.
 \end{equation}
 Clearly, we have 
 \begin{equation}
 \label{err_ineq}
 \terr_\ell\le M(t)\gl_0(t)+\tau_0 \sum_{k\ne \ell}A_k(t).
 \end{equation}
 

For a fixed $t$, and a positive number  $\wt \ep_1$, 
we let $\cG_t$ and $\cG_{t, k}$ denote the sets defined by 
\begin{equation}
\label{def_cGk}
\cG_t:=\{\eta: \; |\wt V_x(t, \eta)|>\wt \ep_1\}, \; \cG_{t, k}:=\{\eta \in \cG_t: \; |\eta-\phi_k^\gp(t)|< \frac {\ga}{\gs(t)}\}.  
\end{equation}
Note that $\cG_t$ and $\cG_{t, k}$ depend on $\wt \ep_1$, and for simplicity of presentation, we drop  
$\wt \ep_1$ from them. Also observe that $\cG_{t, k}=\cG_t\cap \{\eta: \; (t, \eta)\in Z_ k\}$.  
Denote 
\begin{equation}
\label{def_max_eta}
\wh \eta_0:=0, \; \wh \eta_\ell =\wh \eta_\ell(t):={\rm argmax}_{\eta \in\mathcal{G}_{t, \ell}  }|\wt V_x(t,\eta)|, ~~ \ell=1,\cdots, K.
\end{equation}

\begin{theo}
\label{theo:main_1st} Let $x(t)\in {\cB}_{\vep_1, \vep_2}$ for some $\vep_1, \vep_2>0$ and $g$ be a window function.  Suppose  $\gs(t)>0$ 
satisfies \eqref{separated_cond_1st} and  
that 
\begin{equation}
\label{theo1_cond1}
2M(t)\big(\tau_0+\gl_0(t)\big)\le \mu(t). 
\end{equation} 
Let $\wt \ep_1=\wt \ep_1(t)>0$ be a function satisfying   
\begin{equation}
\label{cond_ep1}
M(t)\big(\tau_0+\gl_0(t)\big)\le 
\wt \ep_1 \le \mu(t)-M(t)\big(\tau_0+\gl_0(t)\big).
\end{equation}
Then the following statements hold.
\begin{enumerate}
\item[{\rm (a)}] Let $\cG_t$ and $\cG_{t, k}$ be the sets defined by \eqref{def_cGk} for some $\wt \ep_1$ satisfying \eqref{cond_ep1}.   Then $\cG_t$ can be expressed as a disjoint union of exactly $K+1$ non-empty sets $\cG_{t, k}, 0\le k\le K$. 

\item[{\rm (b)}] Let $\wh \eta_\ell$ be defined by \eqref{def_max_eta}.  Then for $\ell=1, \cdots, K$, 
\begin{equation}\label{phi_est}
|\wh{\eta}_{\ell}(t)-\phi_{\ell}^{'}(t)|\le \bd_{1, \ell}:=\frac{1}{\gs(t)} |\wh g|^{-1}\big(1-\frac {2 \; \terr_\ell(t)}{A_\ell(t)}\big), 
\end{equation}
where $\terr_\ell$ is defined by \eqref{def_err}. 
\item[{\rm (c)}]  For $\ell=0, 1, \cdots, K$,
\begin{equation}
\label{comp_xk_est}
\big|\wt V_{x}(t, \wh \eta_\ell)-x_\ell(t)\big|
\le\bd_{2, \ell}:= \terr_\ell(t)+2\pi I_1 A_\ell(t) |\wh g|^{-1}\big(1-\frac {2 \; \terr_\ell(t)}{A_\ell(t)}\big). 
\end{equation}
\item[{\rm (d)}] 
 If in addition the window function $g(t)\ge 0$ for $t\in \R$, then for $\ell=0, 1, \cdots, K$,
\begin{equation}
\label{abs_IA_est}
\big| |\wt V_x(t,\widehat{\eta}_{\ell})|-A_{\ell}(t) \big|\le \terr_\ell(t). 
\end{equation} 
\end{enumerate}
\end{theo}

Note \eqref{theo1_cond1} holds if $\vep_1, \vep_2, \tau_0$ are small enough. \eqref{abs_IA_est} gives an error bound for IA and trend estimate
 $$
A_{\ell}(t) \approx |\wt V_x(t,\widehat{\eta}_{\ell})|.
$$
The proof of Theorem \ref{theo:main_1st} will delayed to the end of this section. 

\begin{mrem}\label{rem1}
Here we remark that $\terr_\ell(t)< \frac 12 A_\ell(t)$ if $2M(t)\big(\tau_0+\gl_0(t)\big)\le  \mu(t)$, and hence, $|\wh g|^{-1}\big(1-\frac {2 \; \terr_\ell(t)}{A_\ell(t)}\big)$ is well defined. Indeed, from \eqref{err_ineq}, 
\begin{eqnarray*}
 \terr_\ell(t)\hskip -0.6cm&& \le M(t)\gl_0(t)+\tau_0 \sum_{k\ne \ell}A_k(t)
\le \frac 12 \mu(t)-M(t)\tau_0+\tau_0 \sum_{k\ne \ell}A_k(t)\\
&& < \frac 1 2 A_\ell(t)-M(t)\tau_0+\tau_0 \sum_{k=1}^KA_k(t)= \frac 1 2 A_\ell(t). 
\end{eqnarray*}
\hfill $\blacksquare$ 
\end{mrem}

When $\wh g(\xi)$ is supported in $[-\ga, \ga]$, we can set $\tau_0$ in Theorem \ref{theo:main_1st} to be zero.  Thus  the condition in \eqref{theo1_cond1} is reduced to  $2M(t)\lambda_0(t)\le \mu(t)$. In addition, the error $\terr_\ell(t)$ in \eqref{def_err} is simply $M(t)\gl_0(t)$.  To summarize, we have the following corollary.  

\begin{cor}
\label{cor:main_1st} Let $x(t)\in {\cB}_{\vep_1, \vep_2}$ for some $\vep_1, \vep_2>0$ and 
$g$ be a window function with supp($\wh g)\subseteq [-\ga, \ga]$. 
Suppose $\gs(t)>0$ satisfies \eqref{separated_cond_1st}. 
If $\vep_1, \vep_2$ are small enough such that $2M(t)\gl_0(t)\le \mu(t)$, then we have the following.
\begin{enumerate}
\item[{\rm (a)}]  Let $\cG_t$ and $\cG_{t, k}$ be the sets defined by \eqref{def_cGk} for some $\wt \ep_1>0$. 
If  $\wt \ep_1$ satisfies  $M(t)\gl_0(t)\le 
\wt \ep_1 \le \mu(t)-M(t)\gl_0(t)$,   then $\cG_t$ can be expressed as a disjoint union of exactly $K+1$ non-empty sets $\cG_{t, k}, 0\le k\le K$. 

\item[{\rm (b)}] Let $\wh \eta_\ell$ be defined by \eqref{def_max_eta}.  Then for $\ell=1, 2, \cdots, K$,
\begin{equation}\label{cor_phi_est}
|\wh{\eta}_{\ell}(t)-\phi_{\ell}^{'}(t)|\le \frac{1}{\gs(t)} |\wh g|^{-1}\big(1-\frac {2 \; M(t)\gl_0(t)}{A_\ell(t)}\big).  
\end{equation}
\item[{\rm (c)}] For $\ell=0, 1, \cdots, K$,
\begin{equation}
\label{cor_comp_xk_est}
\big|\wt V_{x}(t, \wh \eta_\ell)-x_\ell(t)\big|
\le M(t)\gl_0(t)+2\pi I_1 A_\ell(t) |\wh g|^{-1}\big(1-\frac {2 \; M(t)\gl_0(t)}{A_\ell(t)}\big). 
\end{equation}
\item[{\rm (d)}] 
 If in addition the window function $g(t)\ge 0$ for $t\in \R$, then for $\ell=0, 1, \cdots, K$,
\begin{equation}
\label{cor_abs_IA_est}
\big| |\wt V_x(t,\widehat{\eta}_{\ell})|-A_{\ell}(t) \big|\le M(t)\gl_0(t). 
\end{equation} 
\end{enumerate}
\end{cor}

\bigskip 


Next let us consider the case that $g(t)$ is the Gaussian window function given by  \eqref{def_g}. In this case  $\wh g(\xi)=e^{-2\pi^2 \xi^2}$. 
Thus 
$\terr_\ell(t)$ defined by \eqref{def_err}is 
\begin{equation}
 \label{def_err_gaussian}
 \terr_\ell(t)=\vep_1  I_1\gs (t)M(t) +\pi  \vep_2 I_2 \gs^2 (t)M(t)+\sum_{k\ne \ell}A_k(t)e^{-2\pi^2 \big(\ga (2|\ell-k|-1)\big)^2}.
 \end{equation}

For this $g$, we have 
$$
|\wh g|^{-1}(\xi)=\wh g(\xi)^{-1}= \frac1{\pi \sqrt 2} \sqrt{-\ln \xi}, \;  0<\xi <1. 
$$
Hence the error bound $\bd_{1, \ell}$ in \eqref{phi_est} is 
$$
\bd_{1, \ell}=\frac{1}{\gs(t)} |\wh g|^{-1}\big(1-\frac {2 \; \terr_\ell(t)}{A_\ell(t)}\big)=
\frac{1}{\gs(t)\pi \sqrt 2} \sqrt{-\ln \big(1-\frac {2 \; \terr_\ell(t)}{A_\ell(t)}\big)}. 
$$
Using the fact $-\ln(1-t)< 2t$ for $0<t<\frac12$ and assuming $\vep_1, \vep_2, \tau_0$ are small enough such that $\frac {2 \; \terr_\ell(t)}{A_\ell(t)}<\frac12$, we have 
\begin{eqnarray*}
\bd_{1, \ell}\hskip -0.6cm &&< \frac{\sqrt 2}{\gs(t)\pi} \sqrt{\frac{\terr_\ell(t)}{A_\ell(t)}}. 
\end{eqnarray*}
In this case the error bound $\bd_{2, \ell}$ in \eqref{comp_xk_est} for component recovery satisfies 
$$
\bd_{2, \ell}<\terr_\ell(t)+2\sqrt 2I_1 \sqrt{A_\ell(t)\terr_\ell(t)}. 
$$


To summarize, we have the following corollary. 
\begin{cor}
\label{cor:sinu_cor2} 
Let $x(t)\in {\cB}_{\vep_1, \vep_2}$ for some $\vep_1, \vep_2>0$. Suppose the conditions in 
Theorem \ref{theo:main_1st} for the Gaussian window function $g$ given by \eqref{def_g} are satisfied and that $\terr_\ell(t)< \frac 14 A_\ell(t)$. Then Part {\rm (a)} of  Theorem \ref{theo:main_1st} holds; and  with $\wh \eta_\ell$ defined by \eqref{def_max_eta}, we have  
\begin{eqnarray}
\label{phi_est_gaussian}
&&|\wh{\eta}_{\ell}(t)-\phi_{\ell}^{'}(t)|< 
\frac{\sqrt 2}{\gs(t)\pi} \sqrt{\frac{\terr_\ell(t)}{A_\ell(t)}}, \quad  \ell=1, 2, \cdots, K, \\
\label{comp_xk_est_gaussian}
&& \big|\wt V_{x}(t, \wh \eta_\ell)-x_\ell(t)\big|
< \terr_\ell(t)+2\sqrt 2I_1 \sqrt{A_\ell(t)\terr_\ell(t)}, \quad \ell=0, 1, \cdots, K, \\
\label{cor_abs_IA_est_gaussian}
&&\big| |\wt V_x(t,\widehat{\eta}_{\ell})|-A_{\ell}(t) \big|\le \terr_\ell(t), \quad \ell=0, 1, \cdots, K,
\end{eqnarray}
where $\terr_\ell(t)$ is defined by \eqref{def_err_gaussian}. 
\end{cor}

\begin{mrem} Observe that a smaller  $\gs(t)$ results in smaller error bounds in \eqref{comp_xk_est_gaussian} and \eqref{cor_abs_IA_est_gaussian} for component recovery and instantaneous amplitude recovery respectively.  Note that the IF estimate error bound in  \eqref{phi_est_gaussian} is 
\begin{eqnarray*}
 &&\frac{\sqrt 2}{\gs(t)\pi}\sqrt{\frac{\terr_\ell(t)}{A_\ell(t)}}
=\frac{\sqrt 2}{\pi} \Big(\frac{M(t)}{A_\ell(t)\gs(t)}\vep_1  I_1 +\frac{M(t)}{A_\ell(t)}\pi  \vep_2 I_2 +\frac1 {\gs^2(t)A_\ell(t)}\sum_{k\ne \ell}A_k(t) 
e^{-2\pi^2 \big(\ga (2|\ell-k|-1)\big)^2}\Big)^{1/2}. 
\end{eqnarray*}
Thus for IF estimate, a larger $\gs(t)$ actually gives a smaller error bound. Of course $\gs(t)$ cannot be arbitrarily large due to the restriction in  \eqref{theo1_cond1}. 
 \end{mrem}

\bigskip 

Finally, in this section,  we give the proof of Theorem \ref{theo:main_1st}. Here we consider the case $g$ is non-bandlimited. The proof of Theorem \ref{theo:main_1st} when $g$ is bandlimited, which is Corollary \ref{cor:main_1st}, is the same but simpler. 

\bigskip 

{\bf Proof  of  Theorem \ref{theo:main_1st}(a)}.  Clearly $\cup _{k=0}^K \cG_{t, k}\subseteq \cG_t$. 
Next we show  $\cG_t\subseteq \cup _{k=0}^K \cG_{t, k}$.  Let $\eta \in \cG_t$. 
Assume $\eta\not \in \cup _{k=0}^K \cG_{t, k}$. That is $(t, \eta)\not \in \cup _{k=0}^K Z_k$.  
Then by the definition of $Z_k$ in \eqref{def_Zk}, we have  
$|\wh g\big(\gs(t)(\eta-\phi_k^\gp(t))\big)|\le \tau_0$. Hence, by \eqref{STFT_approx_1st} and \eqref{rem0_est}, we have 
\begin{eqnarray*}
|\wt V_x(t, \eta)|\hskip -0.6cm &&\le \sum_{k=0}^K |x_k(t) \wh g\big(\gs(t)(\eta-\phi_k^\gp(t))\big)| + |\rem_0|
\\
&&\le \tau _0 \sum_{k=0}^K A_k(t)  +M(t)\lambda_0(t)\le \wt \ep_1, 
\end{eqnarray*}
a contradiction to the assumption $|\wt V_x(t, \eta)|>\wt \ep_1$. Thus $(t, \eta)\in Z_\ell$ for some $\ell$. 
This shows that $\eta \in \cG_{t, \ell}$. Hence $\cG_t=\cup _{k=0}^K \cG_{t, k}$. 
Since $Z_k, 0\le k\le K$ are not overlapping, we know  $\cG_{t, k}, 0\le k\le K$ are disjoint.  

To show that $\cG_{t, \ell}$ is non-empty, 
it is enough to show $\phi^\gp_\ell(t)\in \cG_t$ which implies $\phi^\gp_\ell(t)\in \cG_{\ell, t}$ since 
$(t, \phi\rq{}_k(t))\in Z_k$.   
From 
\begin{eqnarray*}
|\wt V_x(t, \phi^\gp_\ell(t))|\hskip -0.6cm &&\ge \sum_{k=0}^K |x_k(t) \wh g\big(\gs(t)(\phi^\gp_\ell(t)-\phi_k^\gp(t))\big)| - |\rem_0|
\\
&&\ge  |x_\ell(t) \wh g(0)| - \sum_{k\ne \ell} |x_k(t)| \tau_0  -M(t)\lambda_0(t)\\
&& \ge  A_\ell(t) - M(t)\tau_0  -M(t)\lambda_0(t)\\
&& \ge  \mu(t) - M(t)\big(\tau_0 +\lambda_0(t)\big) \ge \wt \ep_1, 
\end{eqnarray*}
we conclude that $\phi^\gp_\ell(t)\in \cG_t$. Therefore, the statements in (a) hold. 
\hfill $\blacksquare$

\bigskip 

{\bf Proof  of  Theorem \ref{theo:main_1st}(b)}. For any $\eta\in \cG_{t, \ell}$, from \eqref{STFT_approx_1st}, we have 
\begin{eqnarray*}
\nonumber && \big|\wt V_{x}(t, \eta)-x_\ell(t) \wh g\big(\gs(t)(\eta-\phi_\ell^\gp(t) \big)|\\
&&=\Big|\sum_{k\ne \ell} x_k(t)  \wh g\big(\gs(t)(\eta-\phi_k^\gp(t) \big) +\rem_0\Big|\\
&&\le  \sum_{k\ne \ell} A_k(t) |\wh g\big(\ga (2|\ell-k|-1\big)| + M(t)\gl_0(t), 
\end{eqnarray*}
where the last inequality follows from \eqref{rem0_est}, \eqref{zone_est} and the assumption that $|\wh g(\xi)|$ is decreasing for $\xi\ge 0$. 
Thus 
\begin{equation}
\label{est_xk_STFT}
 \big|\wt V_{x}(t, \eta)-x_\ell(t) \wh g\big(\gs(t)(\eta-\phi_\ell^\gp(t) \big)|\le \terr_\ell(t). 
\end{equation}
Hence, letting $\eta=\phi_\ell^\gp(t)$, we have 
\begin{equation}
\label{est_xk_STFT1}
\big|\wt V_{x}(t, \phi_\ell^\gp(t))\big|\ge |x_\ell(t) \wh g(0)|-\terr_\ell(t) =A_\ell(t) - \terr_\ell(t), 
\end{equation}
since $\wh g(0)=1$. By the definition of $\wh \eta_\ell$ and \eqref{est_xk_STFT} again, we have 
$$
\big|\wt V_{x}(t, \phi_\ell^\gp(t))\big| \le \big|\wt V_{x}(t, \wh \eta_\ell)\big| \le 
\big|x_\ell(t) \wh g\big(\gs(t)(\wh \eta_\ell-\phi_\ell^\gp(t) \big)\big|+ \terr_\ell(t)
$$
This, together with \eqref{est_xk_STFT1}, implies 
$$
A_\ell(t) - \terr_\ell(t)\le A_\ell(t) \big|\wh g\big(\gs(t)(\wh \eta_\ell-\phi_\ell^\gp(t) \big)\big|+ \terr_\ell(t), 
$$
or equivalently
$$
0<1- \frac{2\; \terr_\ell(t)}{A_\ell(t)}\le  \big|\wh g\big(\gs(t)(\wh \eta_\ell-\phi_\ell^\gp(t) \big)\big|.  
$$
Then \eqref{phi_est} follows from the above inequality and that $|\wh g(\xi)|$ is decreasing on $(0, \infty)$. 
\hfill $\blacksquare$

\bigskip 

{\bf Proof  of  Theorem \ref{theo:main_1st}(c)}.
From \eqref{est_xk_STFT}, we have 
\begin{eqnarray*}
&&\big|\wt V_{x}(t, \wh \eta_\ell)-x_\ell(t)\big|\\
&&\le \big|\wt V_{x}(t, \wh \eta_\ell)-x_\ell(t) \wh g\big(\gs(t)(\wh \eta_\ell-\phi_\ell^\gp(t) \big)\big|
+\big |x_\ell(t) \wh g\big(\gs(t)(\wh \eta_\ell-\phi_\ell^\gp(t) \big)-x_\ell(t)\big|
\\
&&\le \terr_\ell(t)+A_\ell(t) \Big| \int_\R g(\tau)\Big(e^{-i2\pi \gs(t)(\wh \eta_\ell-\phi_\ell^\gp(t) )\tau }-1\Big) d\tau \Big|\\
&&\le \terr_\ell(t)+A_\ell(t) \int_\R |g(\tau)| \; 2\pi \gs(t)\big|\wh \eta_\ell-\phi_\ell^\gp(t)\big | |\tau| d\tau\\
&&\le \terr_\ell(t)+A_\ell(t) 2\pi |\wh g|^{-1}\big(1-\frac {2 \; \terr_\ell(t)}{A_\ell(t)}\big) \int_\R |g(\tau)|  |\tau| d\tau\\
&&\le  \terr_\ell(t)+2\pi I_1 A_\ell(t) |\wh g|^{-1}\big(1-\frac {2 \; \terr_\ell(t)}{A_\ell(t)}\big). 
\end{eqnarray*}
This shows \eqref{comp_xk_est}. 
\hfill $\blacksquare$
\bigskip 

{\bf Proof  of  Theorem \ref{theo:main_1st}(d)}. By \eqref{est_xk_STFT1}, we have 
\begin{equation}
\label{est_xk_IA1}
\big|\wt V_{x}(t, \wh \eta_\ell)\big|\ge 
\big|\wt V_{x}(t, \phi_\ell^\gp(t))\big|\ge A_\ell(t) - \terr_\ell(t). 
\end{equation}
On the other hand, when $g(t)\ge 0$, we have $|\wh g(\xi)|\le 1$ for any $\xi\in \R$. This fact and 
\eqref{est_xk_STFT} imply 
 \begin{equation}
\label{est_xk_IA2}
\big|\wt V_{x}(t, \wh \eta_\ell)\big| \le 
\big|x_\ell(t) \wh g\big(\gs(t)(\wh \eta_\ell-\phi_\ell^\gp(t) \big)\big|+ \terr_\ell(t)\le A_\ell(t) +\terr_\ell(t).
\end{equation}
\eqref{abs_IA_est} follows from \eqref{est_xk_IA1} and \eqref{est_xk_IA2}. This completes the proof  of  Theorem \ref{theo:main_1st}(d). 
\hfill $\blacksquare$

\section{Linear chirp-based model}

We consider multicomponent signals $x(t)$ of \eqref{AHM} with $A_k(t)$  satisfying \eqref{condition1} and 
\begin{eqnarray}
\label{cond_basic_2nd} && 
  \phi_k(t)\in C^3(\R),  
  \phi^{\gp\gp}_k(t) \in L_\infty(\R),\\
\label{condition2}
&& |\phi_k^{(3)}(t)|\le \vep_3, \; t\in \R, \; 1\le k\le K, 
\end{eqnarray}
for some positive number $\vep_3$.  As shown below, conditions \eqref{condition1} and \eqref{condition2} imply that 
each component $x_k(t)$ is well approximated locally by linear chirps.  
We say $s(t)$ is  a linear chirp (also called a linear frequency modulated (LFM) signal) if  
\begin{equation}
\label{def_chip_At}
s(t)=A e^{i2\pi \phi(t)}=A e^{i2\pi (ct +\frac 12 r t^2)}, 
\end{equation}
for some constants $c, r$ with $r\not=0$. 
\bigskip

For a given $t$, we use $G_k(\xi)$ to denote the Fourier transform of $e^{i\pi \gs(t)\phi^{\gp\gp}_k(t) \tau^2}g(\tau), \tau\in \R$, namely, 
\begin{equation}
\label{def_Gk}
G_k(\xi)=G_{k, t}(\xi):=
\int_{\R} e^{i\pi\gs^2(t) \phi^{\gp\gp}_k(t) \tau^2}g(\tau) e^{-i2\pi \xi \tau}d\tau.   
\end{equation} 
Note that $G_k(\xi)$ depends on $t$ also if $\phi^{\gp\gp}_k(t)\not=0$. We drop $t$ in $G_{k, t}$ 
for simplicity. For a window function $g$, 
define 
\begin{equation}
\label{def_PFT}
\wb g(\xi, \gl):={\cal F}\Big(e^{-i\pi \gl \tau^2}g(\tau)\Big)\big(\xi)=
\int_{\R} g(\tau) e^{-i2\pi \xi \tau-i\pi \gl \tau^2}d\tau.    
\end{equation}

\begin{mdef} {\rm (}{\bf Admissible window function}{\rm )} \; 
A window function $g\in L_1(\R)\cap L_2(\R)$ is called an admissible window function if for any $\gl\in \R$, $|\wb g(\xi, \gl)|$ is even in $\xi$ and is decreasing for $\xi\ge 0$. 
\end{mdef}

If $g$ is an admissible window function, then for  $G_k(\xi)$ defined by \eqref{def_Gk},  $|G_k(\xi)|$ is even and decreasing for $\xi\ge 0$. We use  $|G_k|^{-1}(\xi)$ to denote the inverse function of   $|G_k(\xi)|, \xi \ge 0$. 
In this section, we assume window function $g$ is admissible. 

\bigskip 

For each component $x_k(t)=A_k(t)e^{i2\pi \phi_k(t)}$ in \eqref{AHM}, $0\le k \le K$, we write $x_k(t+\tau)$ as 
\begin{eqnarray*}
x_k(t+\tau)\hskip -0.6cm &&=x_k(t)e^{i2\pi (\phi^\gp_k(t)\tau+\frac 12\phi^{\gp \gp}_k(t)\tau^2) }
+x_{k, {\rm rl}}(t, \tau),
\end{eqnarray*}
where $x_{k, {\rm rl}}(t, \tau)=x_k(t+\tau) -x_k(t)e^{i2\pi (\phi^\gp_k(t)\tau+\frac 12\phi^{\gp \gp}_k(t)\tau^2) }$. Note that for a fixed $t$, $x_k(t)e^{i2\pi (\phi^\gp_k(t)\tau+\frac 12\phi^{\gp \gp}_k(t)\tau^2)}$  is a linear chirp in $\tau$ (assuming $\phi^{\gp \gp}_k(t)\not=0$). We use $x_k(t)e^{i2\pi (\phi^\gp_k(t)\tau+\frac 12\phi^{\gp \gp}_k(t)\tau^2) }$ to approximate $x(t+\tau)$ at a local time $t$ in the following way: 
\begin{eqnarray}
\nonumber
\wt V_{x}(t, \eta)\hskip -0.6cm &&=\sum_{k=0}^K \int_\R x_k(t+\tau)\frac 1{\gs(t)}g(\frac \tau{\gs(t)}) e^{-i2\pi \eta \tau}d\tau\\
\label{STFT_approx0}
&&=\sum_{k=0}^K \int_\R x_k(t)e^{i2\pi (\phi_k^\gp(t) \tau +\frac12\phi^{\gp\gp}_k(t) \tau^2)}\frac 1{\gs(t)}g(\frac \tau{\gs(t)}) e^{-i2\pi \eta \tau}d\tau +\err_0\\
\label{STFT_approx}
&&=\sum_{k=0}^K x_k(t)  
G_k\big(\gs(t)(\eta-\phi_k^\gp(t) )\big)+\err_0, 
\end{eqnarray}
where 
\begin{eqnarray}
\label{def_err0}
&&\err_0:=\sum_{k=0}^K \int_\R x_{k, {\rm rl}}(t, \tau) \frac 1{\gs(t)}g(\frac \tau{\gs(t)}) e^{-i2\pi \eta \tau}d\tau. 
\end{eqnarray}

Write $x_{k, {\rm rl}}(t, \tau)$ as 
\begin{eqnarray*}
x_{k, {\rm rl}}(t, \tau)\hskip -0.6cm &&:=(A_k(t+\tau)-A_k(t))e^{i2\pi \phi_k(t+\tau)}\\
&&\quad +x_k(t)e^{i2\pi (\phi^\gp_k(t)\tau+\frac 12\phi^{\gp \gp}_k(t)\tau^2) }
\big(
e^{i2\pi (\phi_k(t+\tau)-\phi_k(t)-\phi_k^\gp(t) \tau- \frac 12\phi^{\gp \gp}_k(t)\tau^2)}-1\big).
\end{eqnarray*}
Then by\eqref{condition1} and the fact: 
$$
|e^{i2\pi (\phi_k(t+\tau)-\phi_k(t)-\phi_k^\gp(t) \tau- \frac 12\phi^{\gp \gp}_k(t)\tau^2)}-1|\le 
2\pi \frac 16 \sup_{\eta\in \R}|\phi^{(3)}_k(\eta) \tau^3|
\le \frac \pi 3 \vep_3 |\tau|^3, 
$$
we have  
\begin{eqnarray}
\nonumber |\err_0|\hskip -0.6cm &&\le   \sum_{k=0}^K  \int_\R  \vep_1 A_k(t) | \tau| \frac 1{\gs(t)}|g(\frac \tau{\gs(t)})| d\tau+
\sum_{k=0}^K A_k(t)\int_\R  \frac \pi3 \vep_3 | \tau|^3 \frac 1{\gs(t)}|g(\frac \tau{\gs(t)})|
 d\tau \\
 \label{err0_est} &&=\vep_1 I_1 \gs(t) M(t)  +\frac \pi 3 \vep_3 I_3 \gs^3(t) M(t)=:M(t)\Pi_0(t), 
\end{eqnarray}
where $I_n$ is defined in \eqref{def_In}, and  
\begin{equation}
\label{def_Pi0}
\Pi_0(t):=\vep_1  I_1\gs(t) +\frac \pi 3 \vep_3 I_3 \gs^3(t). 
\end{equation}

Thus if $\vep_1, \vep_3$ are small enough, then $|\err_0|$ is small and hence, $G_k\big(\gs(t)(\eta-\phi_k^\gp(t) )\big)$ provides the time-frequency zone for $\wt V_{x_k}(t, \eta)$.  

Let $0<\tau_0<1$ be a given small number as the threshold for zero. 
Let $\xi_k=\xi_k(t)>0$ denote the unique solution of $|G_k(\xi_k)|=\tau_0.$ 
(Recall that in this section we assume that $|G_k(\xi)|$ is even and decreasing for $\xi\ge 0$.) 
Choose $\ga_k=\ga_k(t)\ge \xi_k(t)$. 
Denote 
\begin{equation}
\label{def_Ok2}
O_k=\{(t, \eta): |\eta-\phi_k^\gp(t)|< \frac {\ga_k}{\gs(t)}, t\in \R\}.   
\end{equation}
Then we have 
\begin{equation}
\label{def_Ok}
|G_k\big(\gs(t)(\eta-\phi_k^\gp)\big)|\le \tau_0 \; \hbox{for any $(t, \eta)\not\in O_k$}. 
\end{equation}

Next we consider, as an example, 
the case when  $g$ is the Gaussian function defined by  \eqref{def_g}
One can obtain for this $g$ (see e.g. \cite{Gibson06,LCHJJ18}),  
\begin{equation}
\label{STFT_LinearChip}
G_{k}(\xi)=\frac {1}{\sqrt{1-i2\pi \phi_k^{\gp\gp}(t)\gs^2(t)}}\;
e^{-\frac{2\pi^2 \xi^2}{1+(2\pi \phi_k^{\gp\gp}(t)\gs^2(t))^2} (1+i2\pi \phi_k^{\gp\gp}(t)\gs^2(t))}, 
\end{equation}
where in this paper, $\sqrt{1-ib}, b\in \RR$, denotes the complex root of 
$1-ib$ lying in the same quadrant as $1-ib$. 
Thus 
\begin{equation}
\label{abs_Gk}
|G_k(\xi)|=\frac 1{\big(1+(2\pi \phi_k^{\gp\gp}(t)\gs^2(t))^2\big)^{\frac 14}}\;
e^{-\frac{2\pi^2}{1+(2\pi \phi_k^{\gp\gp}(t)\gs^2(t))^2}\xi^2}. 
\end{equation} 
Therefore, in this case, assuming $\tau_0 (1+  (2\pi \phi^{\gp\gp}_k(t) \gs^2(t))^2)^{\frac 14}\le 1$, 
\begin{equation}
\label{def_gak_Gaussian}
\xi_k=\sqrt{1+ (2\pi \phi^{\gp\gp}_k(t) \gs^2(t))^2} \; \frac 1{2\pi}\sqrt{2\ln (\frac 1{\tau_0})-\frac 12 \ln(1+ 
 (2\pi \phi^{\gp\gp}_k(t) \gs^2(t))^2)}. 
\end{equation}
As in \cite{LCHJJ18}, one may choose 
\begin{equation}
\label{def_gak}
\ga_k=\ga\big(1+  2\pi |\phi^{\gp\gp}_k(t)| \gs^2(t)\big),
 \end{equation}
where $\ga$ is defined by \eqref{def_ga}. Then $\ga_k\ge \xi_k$.  

\bigskip 

In this section we will assume the multicomponent signal $x(t)$ is well-separated, that is there is $\gs(t)$ such that 
\begin{equation}
\label{cond_no_overlapping} 
O_k\cap O_{\ell}=\emptyset, \quad k\not=\ell. 
\end{equation}
Refer to \cite{LCHJJ18} for the well-separated condition when $g$ is the Gaussian window function. 

From \eqref{cond_no_overlapping} with $\ell=k-1$, 
$$
\phi_k^\gp(t)- \frac {\ga_k}{\gs(t)}\ge \phi_{k-1}^\gp(t)+ \frac {\ga_{k-1}}{\gs(t)}, 
$$
or equivalently 
$$
\gs(t)\big(\phi_k^\gp(t)- \phi_{k-1}^\gp(t)\big)\ge  \ga_k+\ga_{k-1}.  
$$
More generally, we have for $k\ne \ell$, 
$$
\gs(t)\big|\phi_k^\gp(t)- \phi_\ell^\gp(t)\big|\ge  \Upsilon_{k, \ell},  
$$
where 
$$
 \Upsilon_{k, \ell}:=
 \begin{cases}
 \ga_k+\ga_{\ell}+ 2(\ga_{k-1}+\ga_{k-2}+\cdots +\ga_{\ell+1}), \; \hbox{for $k>\ell$}, 
 \\
\ga_k+\ga_{\ell}+ 2(\ga_{\ell-1}+\ga_{\ell-2}+\cdots +\ga_{k+1}), \; \hbox{for $k<\ell$}. 
  \end{cases}
$$
Hence 
for $(t, \eta)\in O_\ell$, the following holds 
\begin{equation}
\label{Ok_phi_ell}
\gs(t)|\eta-\phi_k^\gp(t)|\ge \gs(t)\big|\phi_k^\gp(t) - \phi_\ell^\gp(t)\big|- \gs(t)\big|\eta-\phi_\ell^\gp(t)\big|
\ge \Upsilon_{k, \ell} -\ga_\ell. 
\end{equation}
Denote 
 \begin{eqnarray}
 \label{def_Err}
&& \tErr_\ell=\tErr_\ell(t):=M(t)\Pi_0(t)+\sum_{k\ne \ell}A_k(t) |G_k\big(\Upsilon_{k, \ell} -\ga_\ell\big)|, \\
 && g_0=g_0(t):=\min_{0\le k \le K} |G_k(0)|= \min\{1,  |G_k(0)|, 1\le k\le K\},   
 \end{eqnarray}
 where $\Pi_0(t)$ is defined by \eqref{def_Pi0}. 
 
 Since $\Upsilon_{k, \ell} -\ga_\ell\ge \ga_k$ and $|G_k(\ga_k)|\le \tau_0$, we have 
 $$
 \tErr_\ell\le M(t)\Pi_0(t)+\tau_0 \sum_{k\ne \ell}A_k(t).  
 $$


For a fixed $t$, and a positive number  $\wt \ep_1$, 
we let $\cG_t$ denote the set defined by \eqref{def_cGk}, and we define 
\begin{equation}
\label{def_Hk}
\cH_{t, k}:=\{\eta \in \cG_t: \; |\eta-\phi_k^\gp(t)|< \frac {\ga_k}{\gs(t)}\}.  
\end{equation}
Note that 
$\cH_{t, k}$ depends on $\wt \ep_1$, and for simplicity of presentation, we drop 
$\wt \ep_1$ from it. Also observe that $\cH_{t, k}=\cG_t\cap \{\eta: \; (t, \eta)\in O_ k\}$.  
Denote 
\begin{equation}
\label{def_max_eta_2nd}
\wc \eta_0:=0, \; \wc \eta_\ell :=\wc \eta_\ell(t):={\rm argmax}_{\eta \in\mathcal{H}_{t, \ell}  }|\wt V_x(t,\eta)|, ~~ \ell=1,\cdots, K.
\end{equation}

Let ${\cD}_{\vep_1, \vep_3}$ denote the set of multicomponent signals of \eqref{AHM} satisfying 
\eqref{cond_basic}-\eqref{condition1}, 
\eqref{cond_basic_2nd}, \eqref{condition2}.

\begin{theo}
\label{theo:main_2nd} Let $x(t)\in {\cD}_{\vep_1, \vep_3}$ for some small $\vep_1, \vep_3>0$     
and $g$ be an admission window function.  
Suppose for this $x(t)$, there is a function $\gs(t)>0$ such that \eqref{cond_no_overlapping} holds, and  
 $\vep_1, \vep_3, \tau_0$ are small enough such that 
$$
\tErr_\ell<\frac 12{|G_\ell(0)|} A_\ell(t), \; 
2M(t)\big(\tau_0+\Pi_0(t)\big)\le g_0(t)\mu(t). 
$$ 
Let $\wt \ep_1$ be a function satisfying  
\begin{equation}
\label{cond_ep_2nd}
M(t)\big(\tau_0+\Pi_0(t)\big)\le 
\wt \ep_1 \le g_0(t) \mu(t)-M(t)\big(\tau_0+\Pi_0(t)\big).
\end{equation}
Then the following statements hold. 
\begin{enumerate}
\item[{\rm (a)}] Let $\cG_t$ and $\cH_{t, k}$ be the sets defined by \eqref{def_cGk} and \eqref{def_Hk} respectively for some $\wt \ep_1$ satisfying \eqref{cond_ep_2nd}.   Then $\cG_t$ can be expressed as a disjoint union of exactly $K+1$ non-empty sets $\cH_{t, k}, 0\le k\le K$. 

\item[{\rm (b)}] Let $\wc \eta_\ell$ be defined by \eqref{def_max_eta_2nd}.  Then for $\ell=1, \cdots, K$, 
\begin{equation}\label{phi_est_2nd}
|\wc\eta_{\ell}(t)-\phi_{\ell}^{'}(t)|\le \Bd_{1, \ell}:=\frac{1}{\gs(t)} |G_\ell|^{-1}\big(|G_\ell(0)|-\frac {2 \; \tErr_\ell(t)}{A_\ell(t)}\big), 
\end{equation}
where $\tErr_\ell$ is defined by \eqref{def_Err}. 
\item[{\rm (c)}]  For $\ell=0, \cdots, K$,
\begin{equation}
\label{comp_xk_est_2nd}
\big|\wt V_{x}(t, \wc \eta_\ell)-G_\ell(0) x_\ell(t)\big|
\le \Bd_{2, \ell}:=\tErr_\ell(t)+2\pi I_1 A_\ell(t) |\wh G_\ell|^{-1}\big(|G_\ell(0)|-\frac {2 \; \tErr_\ell(t)}{A_\ell(t)}\big). 
\end{equation}
\end{enumerate}
\end{theo}

The proof of Theorem \ref{theo:main_2nd} will be deferred to the end of this section. 

From Theorem \ref{theo:main_2nd}, we have that $\wc \eta_\ell$ defined by \eqref{def_max_eta_2nd} provides an approximation to $\phi_{\ell}^{'}(t)$. In addition, we have 
the recovery formula
\begin{equation}
\label{comp_xk_est_2nd1}
 x_\ell(t) \approx \frac 1{G_\ell(0)}\wt V_{x}(t, \wc \eta_\ell). 
\end{equation}
For a real $x(t)$, the  recovery formula will be  
\begin{equation}
\label{comp_xk_est_2nd1_real}
 x_\ell(t) \approx 2{\rm Re}\Big(\frac 1{G_\ell(0)}\wt V_{x}(t, \wc \eta_\ell)\Big). 
\end{equation}

In some case $\wc \eta_\ell$ is $\wh \eta_\ell$ defined by \eqref{def_max_eta}. Thus in this case, 
the linear chirp-based model does not provide a more accurate IF approximation than the sinusoidal signal-based model. However, 
the component recovery formula derived from linear chirp local approximation has a factor 
 $\frac 1{G_\ell(0)}$ as shown in \eqref{comp_xk_est_2nd1}, which distinguishes the linear chirp-based model  from the sinusoidal signal-based model.

\bigskip 

Recall that when $g(t)$ is the Gaussian window function given by \eqref{def_g}, then  the corresponding $G_k(\xi)$ defined by \eqref{def_Gk} and $|G_k(\xi)|$ is given by \eqref{abs_Gk}. Observe that if we choose $\ga_k$ as the quantity given by \eqref{def_gak}, then 
$$
|G_k\big(\Upsilon_{k, \ell} -\ga_\ell\big)|\le |G_k\big(\ga_k\big)| \le \frac 1{\big(1+(2\pi \phi_k^{\gp\gp}(t)\gs^2(t))^2\big)^{\frac 14}}\;
e^{-2\pi^2\ga ^2}\le e^{-2\pi^2\ga ^2}. 
$$
Hence the terms $|G_k\big(\Upsilon_{k, \ell} -\ga_\ell\big)|$  in $\tErr_\ell$ defined by \eqref{def_Err} are 
very small if $\ga$ is quite large, say $\ga\ge 1$. In addition, we have 
$$
|G_\ell|^{-1}(\xi)= \frac1{\pi \sqrt 2|G_\ell(0)|^2} \Big(-\ln \frac {\xi}{|G_\ell(0)|}\Big)^{1/2}, \;  0<\xi <|G_\ell(0)|,  
$$
where 
$$
|G_\ell(0)|=\big(1+(2\pi \phi_\ell^{\gp\gp}(t)\gs^2(t))^2\big)^{-\frac 14}.
$$ 
Hence the error bound $\Bd_{1, \ell}$ in \eqref{phi_est_2nd} is 
$$
\Bd_{1, \ell}=\frac{1}{\gs(t)} |G_\ell|^{-1}\big(|G_\ell(0)|-\frac {2 \; \tErr_\ell(t)}{A_\ell(t)}\big)=
\frac{1}{\gs(t)\pi \sqrt 2 |G_\ell(0)|^2} \Big(-\ln \big(1-\frac {2 \; \tErr_\ell(t)}{|G_\ell(0)|A_\ell(t)}\Big)^{1/2}. 
$$
Suppose $\vep_1, \vep_3, \tau_0$ are small enough such that $\tErr_\ell(t)\le \frac 14 |G_\ell(0)|{A_\ell(t)}$. Then applying the fact $-\ln(1-t)< 2t$ for $0<t<\frac12$ again,  we have 
\begin{eqnarray*}
\Bd_{1, \ell}\hskip -0.6cm &&< \frac{\sqrt 2}{\gs(t)\pi} \frac 1{|G_\ell(0)|^{5/2}}\sqrt{\frac{\tErr_\ell(t)}{A_\ell(t)}}.
\end{eqnarray*}
Also, in this case, we have 
\begin{eqnarray*}
\Bd_{2, \ell}\hskip -0.6cm &&< \tErr_\ell(t)+\frac{2\sqrt 2 I_2}{|G_\ell(0)|^{5/2}}\sqrt{A_\ell(t) \tErr_\ell(t)}. 
\end{eqnarray*}

To summarize, we have the following corollary. 
\begin{cor}
\label{cor:chirp_cor3} 
Let $x(t)\in {\cD}_{\vep_1, \vep_2}$ for some $\vep_1, \vep_2>0$. Suppose the conditions in 
Theorem \ref{theo:main_2nd} for the Gaussian window function $g$ given by \eqref{def_g} are satisfied and that $\tErr_\ell(t)< \frac 14 |G_\ell(0)| A_\ell(t)$. Then Part {\rm (a)} of  Theorem \ref{theo:main_2nd} holds; and  with 
$\wc \eta_\ell$ defined by \eqref{def_max_eta_2nd}, we have  
\begin{eqnarray}
\label{phi_est_2nd_gaussian}
&&|\phi_{\ell}^{'}(t)-\wc \eta_{\ell}(t)|< 
\frac{\sqrt 2}{\gs(t)\pi} \big(1+(2\pi \phi_\ell^{\gp\gp}(t)\gs^2(t))^2\big)^{\frac 58}\sqrt{\frac{\tErr_\ell(t)}{A_\ell(t)}}, \quad  \ell=1, 2, \cdots, K, \\
\label{comp_xk_est_2nd_gaussian}
&& \big|x_\ell(t)-\sqrt{1-i2\pi \phi_\ell^{\gp\gp}(t)\gs^2(t)} \; \wt V_{x}(t, \wc \eta_\ell)\big|
< \big(1+(2\pi \phi_\ell^{\gp\gp}(t)\gs^2(t))^2\big)^{\frac 14}\; \tErr_\ell(t)\\
\nonumber &&\qquad \qquad \quad  +2\sqrt 2I_1 \big(1+(2\pi \phi_\ell^{\gp\gp}(t)\gs^2(t))^2\big)^{\frac 78}\sqrt{A_\ell(t)\tErr_\ell(t)}, \quad \ell=0, 1, \cdots, K, 
\end{eqnarray}
where $\tErr_\ell(t)$ is defined by \eqref{def_Err}. 
\end{cor}

\begin{mrem} From \eqref{comp_xk_est_2nd_gaussian}, we will use the following  recovery formula for a real $x(t)$:  
\begin{equation}
\label{comp_xk_est_2nd1_real_gaussian} 
x_\ell(t)\approx 2{\rm Re }\Big(\sqrt{1-i2\pi \phi_\ell^{\gp\gp}(t)\gs^2(t)} \; \wt V_{x}(t, \wc \eta_\ell)\Big). 
\end{equation}
\eqref{comp_xk_est_2nd1_real_gaussian} was first derived in \cite{LCJ20} using linear chirp local approximation. Here we provide an error bound as shown in \eqref{comp_xk_est_2nd_gaussian}.
\end{mrem}

\begin{mrem}Observe that a smaller  $\gs(t)$ results in smaller error bounds in \eqref{comp_xk_est_2nd_gaussian} for component recovery.  Note that the IF estimate error bound in  \eqref{phi_est_2nd_gaussian} is 
\begin{eqnarray*}
 &&\frac{\sqrt 2}{\pi \sqrt{A_\ell(t)}} \big(1+(2\pi \phi_\ell^{\gp\gp}(t)\gs^2(t))^2\big)^{\frac 58}
\Big(\frac{M(t)}{\gs(t) }\vep_1  I_1 +M(t)\frac {\pi}3  \vep_3 I_3 \gs (t)+\frac1 {\gs^2(t)}\sum_{k\ne \ell}A_k(t) |G_k\big(\Upsilon_{k, \ell} -\ga_\ell\big)|
 \Big)^{1/2} . 
\end{eqnarray*}
Thus for IF estimate, if $\phi_\ell^{\gp\gp}(t)$ is quite large, then we should choose a smaller $\gs(t)$, which should not be too small to avoid a possible large error caused by  
$ \frac{M(t)}{\gs(t) }\vep_1  I_1$ or/and $\frac1 {\gs^2(t)}\sum_{k\ne \ell}A_k(t) |G_k\big(\Upsilon_{k, \ell} -\ga_\ell\big)|$. 
\end{mrem}
   
 \begin{mrem}
 Unlike Theorem \ref{theo:main_1st} and Corollaries \ref{cor:main_1st} and \ref{cor:sinu_cor2}, there is no restriction of the  amplitude of $\phi^{\gp\gp}_k(t)$ for the source signal $x(t)$ to be analyzed in Theorem \ref{theo:main_2nd} and Corollary \ref{cor:chirp_cor3}. However, as we see from \eqref{phi_est_2nd}  \eqref{comp_xk_est_2nd} and  \eqref{phi_est_2nd_gaussian}  \eqref{comp_xk_est_2nd_gaussian}, $\phi^{\gp\gp}_k(t)$ still effects the IF estimate error and component recovery error.  
 \end{mrem}

\bigskip

Next we provide the proof of Theorem \ref{theo:main_2nd}, which is similar to that of Theorem \ref{theo:main_1st}.

\bigskip 

{\bf Proof  of  Theorem \ref{theo:main_2nd}(a)}.  Clearly $\cup _{k=0}^K \cH_{t, k}\subseteq \cG_t$. 
Next we show  $\cG_t\subseteq \cup _{k=0}^K \cH_{t, k}$.  Let $\eta \in \cG_t$. 
Assume $\eta\not \in \cup _{k=0}^K \cH_{t, k}$. That is $(t, \eta)\not \in \cup _{k=0}^K O_k$.  
From \eqref{def_Ok}, we have 
$|G_k\big(\gs(t)(\eta-\phi_k^\gp(t))\big)|\le \tau_0$ for any $k$. Hence, by \eqref{STFT_approx} and \eqref{err0_est}, we have 
\begin{eqnarray*}
|\wt V_x(t, \eta)|\hskip -0.6cm &&\le \sum_{k=0}^K |x_k(t) G_k\big(\gs(t)(\eta-\phi_k^\gp(t))\big)| + |\err_0|
\\
&&\le \tau _0 \sum_{k=0}^K A_k(t)  +M(t)\Pi_0(t)\le \wt \ep_1, 
\end{eqnarray*}
a contradiction to the assumption $|\wt V_x(t, \eta)|>\wt \ep_1$. Thus $(t, \eta)\in O_\ell$ for some $\ell$, and hence $\eta \in \cH_{t, \ell}$. 
This shows $\cG_t=\cup _{k=0}^K \cH_{t, k}$. 
Since $O_k, 0\le k\le K$ are not overlapping, we know  $\cH_{t, k}, 0\le k\le K$ are disjoint.  

To show $\cH_{t, \ell}$ is non-empty, 
we need only to show $\phi^\gp_\ell(t)\in \cG_t$ which implies $\phi^\gp_\ell(t)\in \cH_{\ell, t}$.   
From 
\begin{eqnarray*}
|\wt V_x(t, \phi^\gp_\ell(t))|\hskip -0.6cm &&\ge \sum_{k=0}^K |x_k(t) G_k\big(\gs(t)(\phi^\gp_\ell(t)-\phi_k^\gp(t))\big)| - |\err_0|
\\
&&\ge  |x_\ell(t) G_\ell(0)| - \sum_{k\ne \ell} |x_k(t)| \tau_0  -M(t)\Pi_0(t)\\
&& \ge  |G_\ell(0)| A_\ell(t) - M(t)\tau_0  -M(t)\Pi_0(t)\\
&& \ge  g_0(t)\mu(t) - M(t)\big(\tau_0 +\Pi_0(t)\big) \ge \wt \ep_1, 
\end{eqnarray*}
we conclude that $\phi^\gp_\ell(t)\in \cG_t$. Therefore, the statements in (a) hold. 
\hfill $\blacksquare$

\bigskip 

{\bf Proof  of  Theorem \ref{theo:main_2nd}(b)}. For any $\eta\in \cH_{t, \ell}$, from \eqref{STFT_approx}, we have 
\begin{eqnarray*}
\nonumber && \big|\wt V_{x}(t, \eta)-x_\ell(t) G_\ell\big(\gs(t)(\eta-\phi_\ell^\gp(t) \big)|\\
&&=\Big|\sum_{k\ne \ell} x_k(t)G_k\big(\gs(t)(\eta-\phi_k^\gp(t) \big) +\err_0\Big|\\
&&\le  \sum_{k\ne \ell} A_k(t) |G_k(\Upsilon_{k, \ell} -\ga_\ell)|+ M(t)\Pi_0(t), 
\end{eqnarray*}
where the last inequality follows from \eqref{err0_est}, \eqref{Ok_phi_ell} and the assumption that $|G_k(\xi)|$ is decreasing for $\xi\ge 0$. 
Thus 
\begin{equation}
\label{est_xk_STFT_2nd}
 \big|\wt V_{x}(t, \eta)- x_\ell(t) G_\ell\big(\gs(t)(\eta-\phi_\ell^\gp(t) \big)|\le \tErr_\ell(t). 
\end{equation}
Hence, letting $\eta=\phi_\ell^\gp(t)$, we have 
\begin{equation}
\label{est_xk_STFT_2nd1}
\big|\wt V_{x}(t, \phi_\ell^\gp(t))\big|\ge |x_\ell(t) G_\ell(0)|-\tErr_\ell(t) =|G_\ell(0)| A_\ell(t) - \tErr_\ell(t).
\end{equation}
By the definition of $\wc \eta_\ell$ and \eqref{est_xk_STFT_2nd} again, we have 
$$
\big|\wt V_{x}(t, \phi_\ell^\gp(t))\big| \le \big|\wt V_{x}(t, \wc \eta_\ell)\big| \le 
\big|x_\ell(t) G_\ell\big(\gs(t)(\wc \eta_\ell-\phi_\ell^\gp(t)) \big)\big|+ \tErr_\ell(t)
$$
This, together with \eqref{est_xk_STFT_2nd1}, implies 
$$
|G_\ell(0)| A_\ell(t) - \tErr_\ell(t)\le A_\ell(t) \big|G_\ell\big(\gs(t)(\wc \eta_\ell-\phi_\ell^\gp(t)) \big)\big|+ \tErr_\ell(t), 
$$
or equivalently
$$
|G_\ell(0)|- \frac{2\; \tErr_\ell(t)}{A_\ell(t)}\le  \big|G_\ell\big(\gs(t)(\wc \eta_\ell-\phi_\ell^\gp(t)) \big)\big|.  
$$
Since $|G_\ell(\xi)|$ is decreasing for $\xi\ge 0$, we have 
$$
\gs(t)\big| \wc \eta_\ell-\phi_\ell^\gp(t) \big|\le |G_\ell|^{-1}\Big(
|G_\ell(0)|- \frac{2\; \tErr_\ell(t)}{A_\ell(t)}\Big). 
$$
Thus shows \eqref{phi_est_2nd}.
\hfill $\blacksquare$

\bigskip 

{\bf Proof  of  Theorem \ref{theo:main_2nd}(c)}.
From \eqref{est_xk_STFT_2nd}, we have 
\begin{eqnarray*}
&&\big|\wt V_{x}(t, \wc \eta_\ell)-G_\ell(0) x_\ell(t)\big|\\
&&\le \big|\wt V_{x}(t, \wc \eta_\ell)- x_\ell(t) G_\ell\big(\gs(t)(\wc \eta_\ell-\phi_\ell^\gp(t)) \big)\big|
+\big |x_\ell(t) G_\ell\big(\gs(t)(\wc \eta_\ell-\phi_\ell^\gp(t)) \big)-G_\ell(0) x_\ell(t)\big|
\\
&&\le \tErr_\ell(t)+A_\ell(t) \Big| \int_\R g(\tau)\Big(e^{i\pi \gs^2(t)\phi_\ell^{\gp\gp}(t) \tau^2} e^{-i2\pi \gs(t)(\wc \eta_\ell-\phi_\ell^\gp(t) )\tau }-e^{i\pi \gs^2(t)\phi_\ell^{\gp\gp}(t) \tau^2}\Big) d\tau \Big|\\
&&\le \tErr_\ell(t)+A_\ell(t) \int_\R |g(\tau)| \; 2\pi \gs(t)\big|\wc \eta_\ell-\phi_\ell^\gp(t)\big | |\tau| d\tau\\
&&\le \tErr_\ell(t)+A_\ell(t) 2\pi |G_\ell|^{-1}\big(|G_\ell(0)|-\frac {2 \; \tErr_\ell(t)}{A_\ell(t)}\big) \int_\R |g(\tau)|  |\tau| d\tau\\
&&=  \tErr_\ell(t)+2\pi I_1 A_\ell(t) |G_\ell|^{-1}\big(|G_\ell(0)|-\frac {2 \; \tErr_\ell(t)}{A_\ell(t)}\big). 
\end{eqnarray*}
This completes the proof of \eqref{comp_xk_est_2nd}. 
\hfill $\blacksquare$
\bigskip 


\section{Experiments}
In this section we provide some experimental results. 
We focus on the comparison between the performances of sinusoidal signal-based and linear chirp-based  models in component recovery. The readers are referred to \cite{LCJ20} for experiments which compare our linear chirp-based model with SST and the 2nd-order SST in component recovery.  In the first two examples we consider monocomponent signals, where we simple let $\gs=\frac 1{16}$. In all our experiments, we use the Gaussian function $g(t)$ defined by \eqref{def_g} as the window function. The signals considered in our experiments are uniformly sampled. Suppose $t_m= m \gt t, m=0,1, \cdots, N-1$ be the sample points. Then finite sequence $\big(\wh \eta_\ell(t_m)\big)_{ 0\le m<N}$ obtained by \eqref{def_max_eta} or  $\big(\wc \eta_\ell(t_m)\big)_{ 0\le m<N}$ in \eqref{def_max_eta_2nd} is an estimate to IF $\phi_\ell^{\gp}(t_m)$. 
In our experiments  $\wh \eta_\ell(t_m) =\wc \eta_\ell(t_m)$. When we use an linear chirp-based model recovery formula such as  \eqref{comp_xk_est_2nd1_real} or 
\eqref{comp_xk_est_2nd1_real_gaussian}, we need to estimate $\phi^{\gp\gp}_\ell(t)$. \cite{LCJ20} provides a method to estimate $\phi^{\gp\gp}_\ell(t)$. Here we use a five-point formula (see, e.g. \cite{BF11_book})  to  
obtain an approximation to $\phi^{\gp\gp}_\ell(t_m)$.
More precisely, we first smooth $\wc \eta_\ell(t_m),  0\le m<N$ by a filter (here we use cubic B-spline filter $\{\frac 1{16}, \frac 4{16} \frac 6{16} \frac 4{16} \frac 1{16}\}$). After that we apply a five-point formula to  the smoothed $\wc \eta_\ell(t_m),  0\le m<N$ to obtain an approximation to   $\big(\phi^{\gp\gp}_\ell(t_m)\big)_{0\le m<N}$,  denoted as  $\big(r_\ell(t_m)\big)_{0\le m<N}$. 
 Finally this approximation  is smoothed by the  cubic B-spline filter and the resulting sequence, denoted as 
 $\big(\wt r_\ell(t_m)\big)_{0\le m<N}$, is  used in \eqref{comp_xk_est_2nd1_real_gaussian} as $\phi^{\gp\gp}_\ell(t)$. In the following examples, in order to avoid the boundary effect, 
 we provide the recovery errors over $[N/8+1:7N/8]$ in either figures or a table below.
  
\begin{figure}[th]
	\centering
	\begin{tabular}{cc}
		\resizebox{3.0in}{2.0in}{\includegraphics{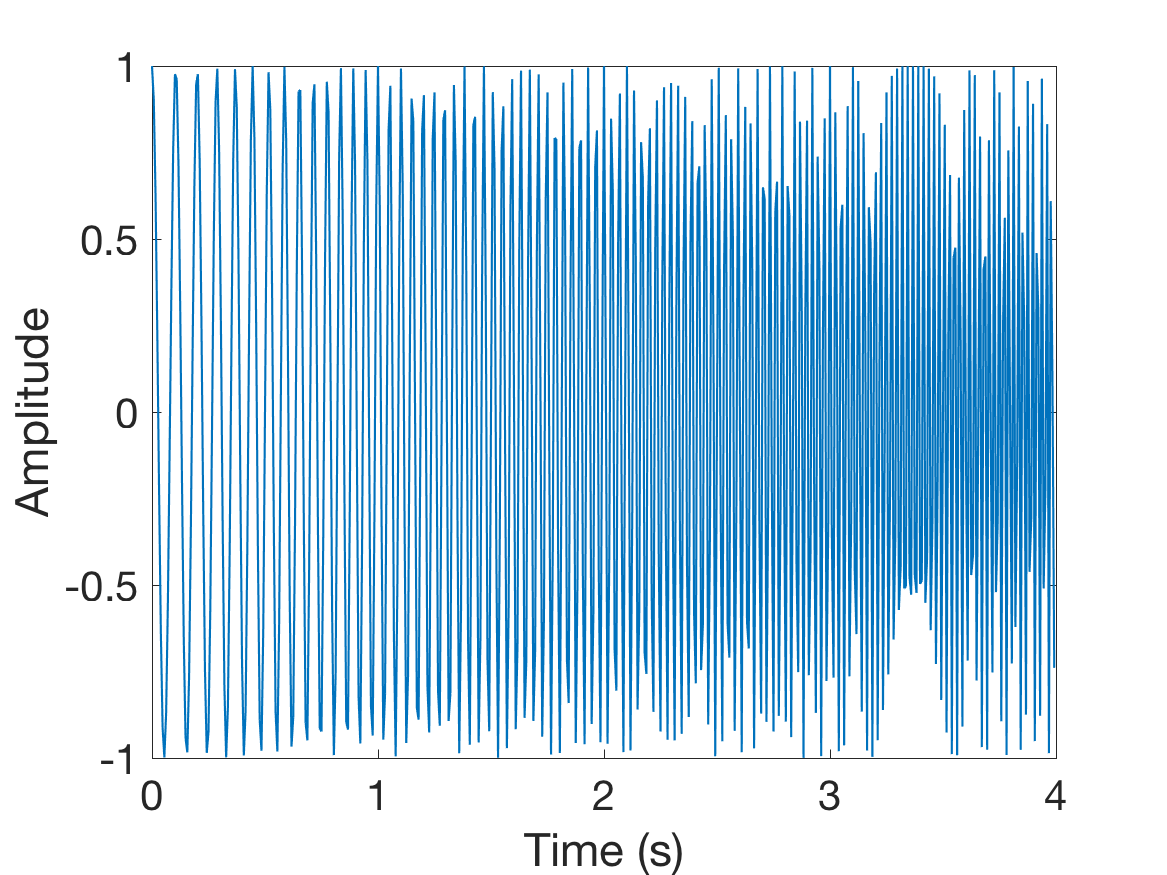}}
		\quad & 
\resizebox{3.0in}{2.0in}{\includegraphics{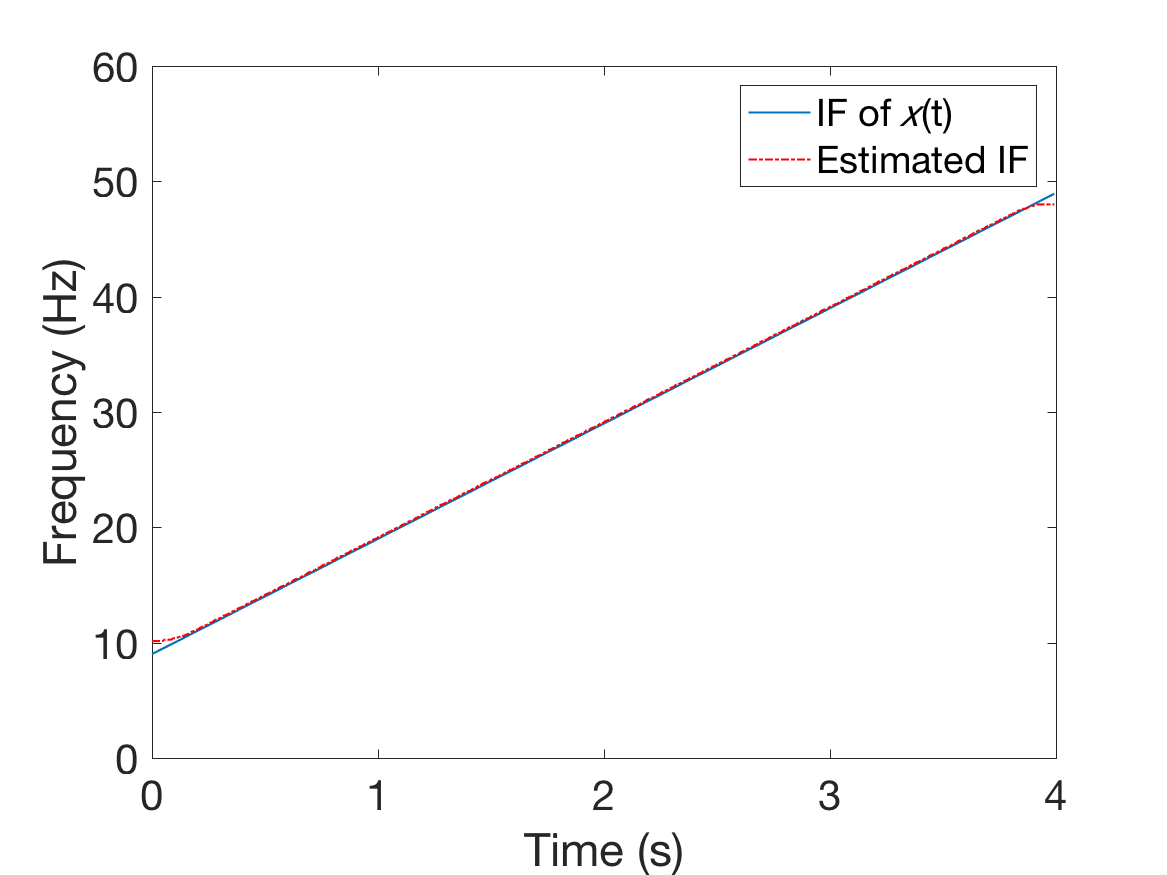}}
\\
 \resizebox{3.0in}{2.0in} {\includegraphics{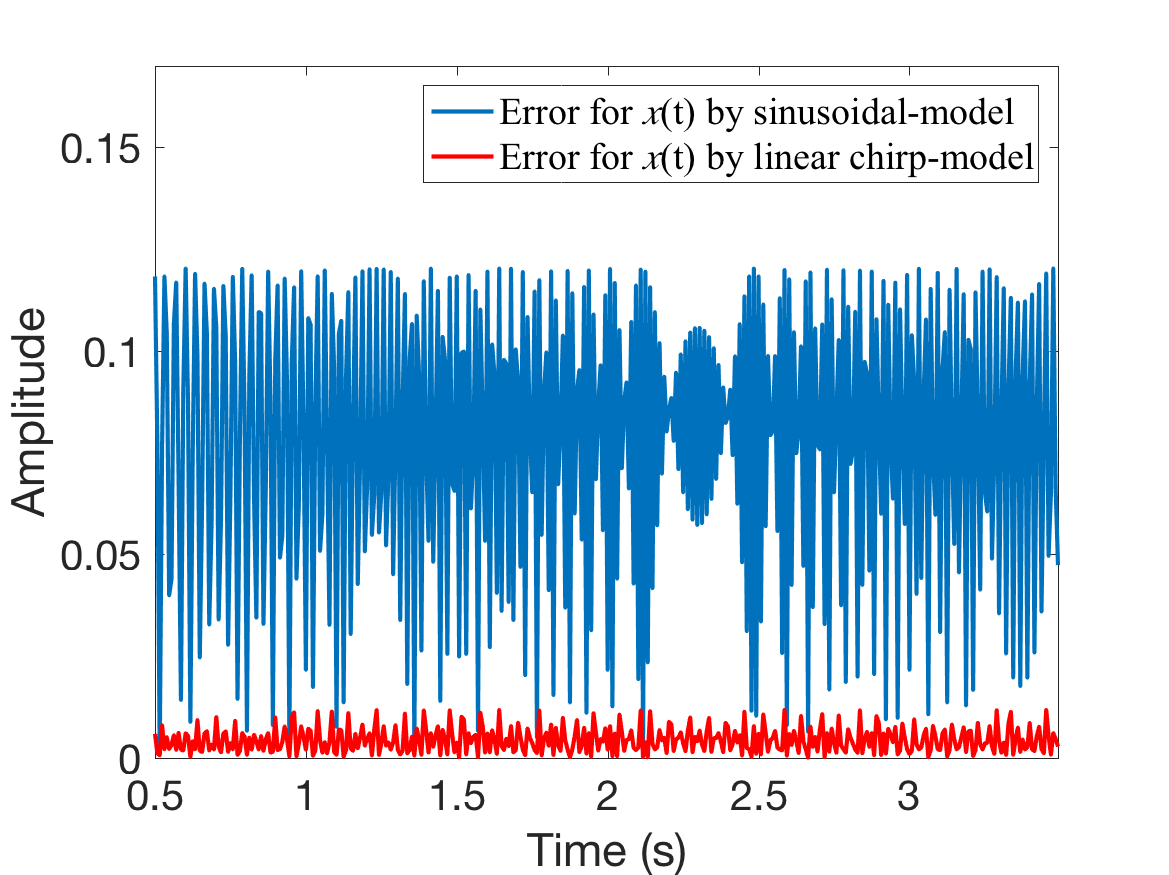}} 
 \quad 
&\resizebox{3.0in}{2.0in} {\includegraphics{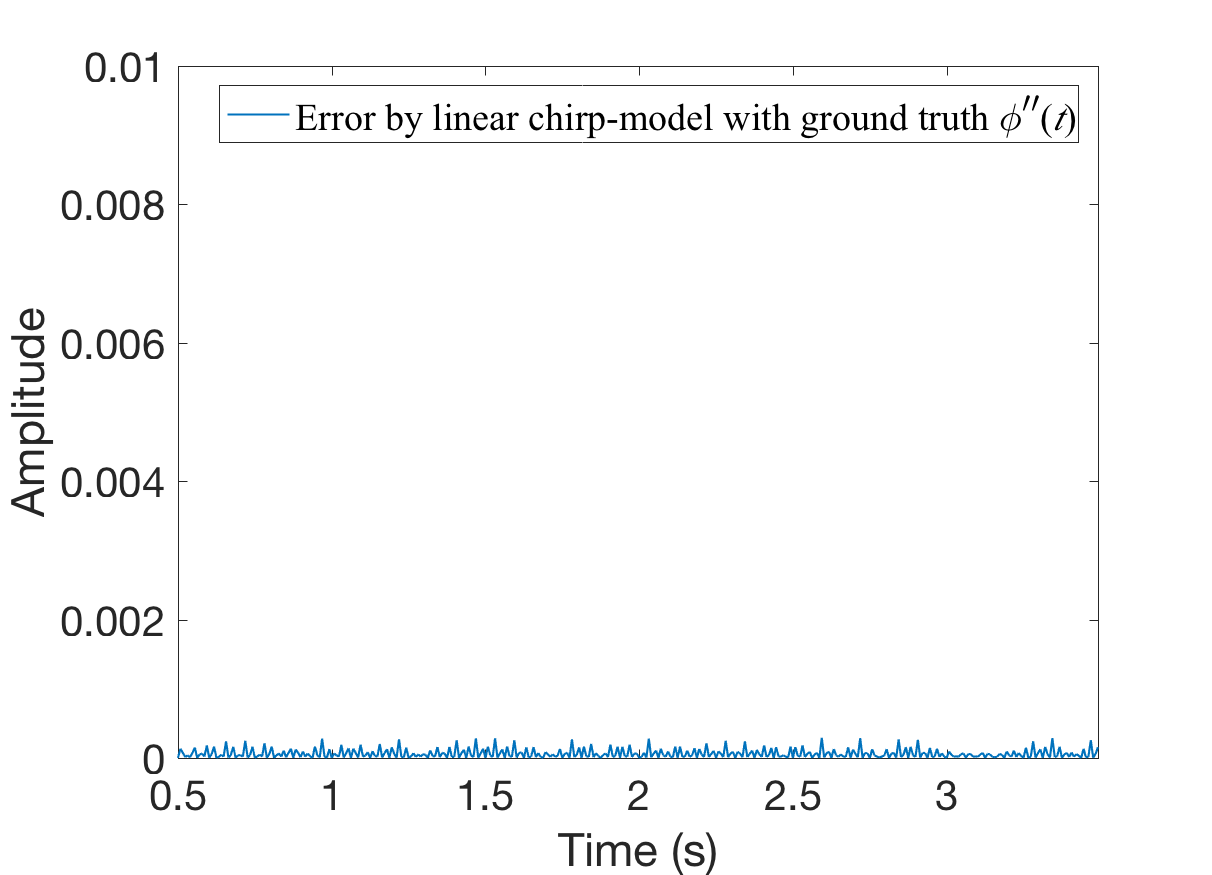}}
	\end{tabular}
	\caption{\small Example of monocomponent signal $x(t)$ in \eqref{one_chirp}. 
		Top-left:  Waveform $x(t)$;  
		Top-right:  IF $\phi\rq{}(t) $ of $x(t)$ (solid line) and estimated IF $\wh \eta(t)$ (red dot-dashed line); 	Bottom-left: Absolute recovery errors for $x(t)$ by  sinusoidal signal-based model (blue line) and by linear chirp-based model with estimated $\phi^{\gp\gp}(t)$ (red line); 
	Bottom-right: Absolute recovery error for $x(t)$ by linear chirp-based model with ground truth $\phi^{\gp\gp}(t)$.}
	\label{fig:one_chirp}
\end{figure}

\begin{figure}[th]
	\centering
	\begin{tabular}{cc}
		\resizebox{3.0in}{2.0in}
{\includegraphics{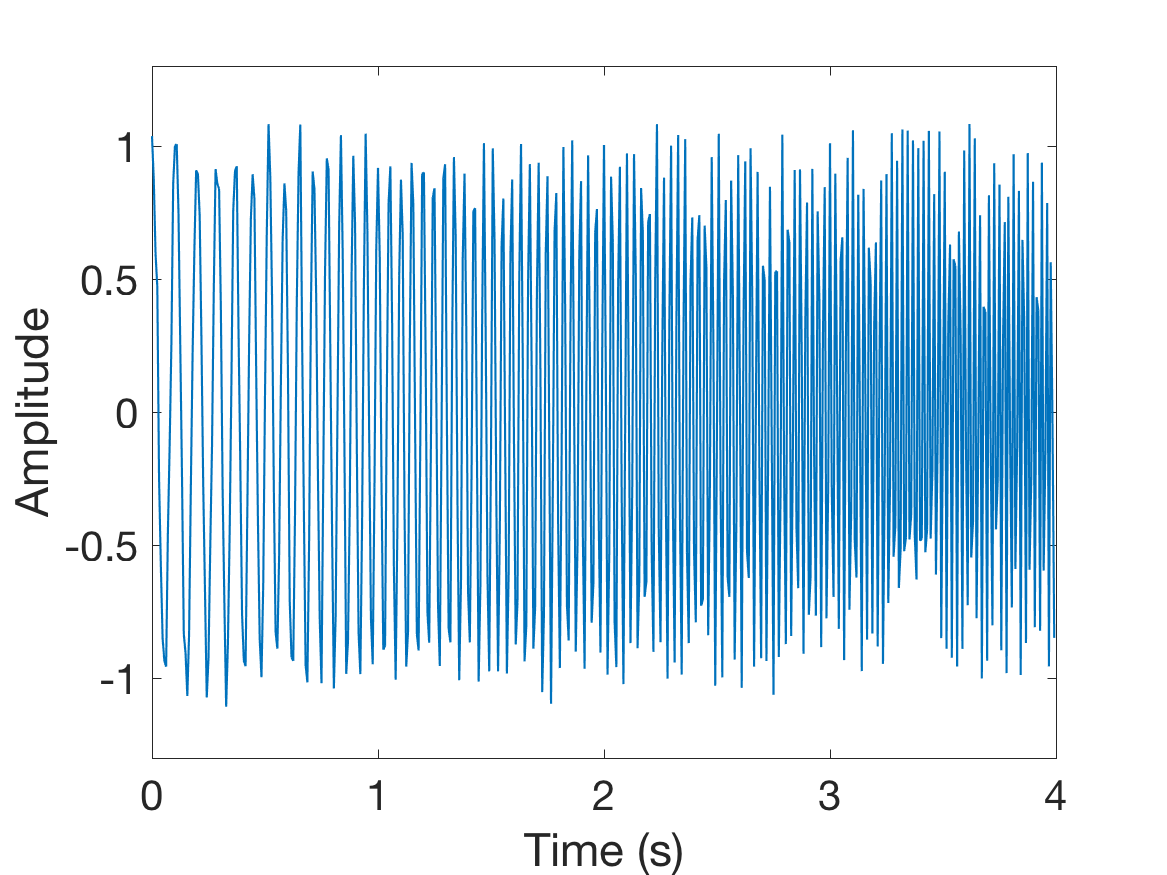}}\quad 
	\quad &\resizebox{3.0in}{2.0in}{\includegraphics{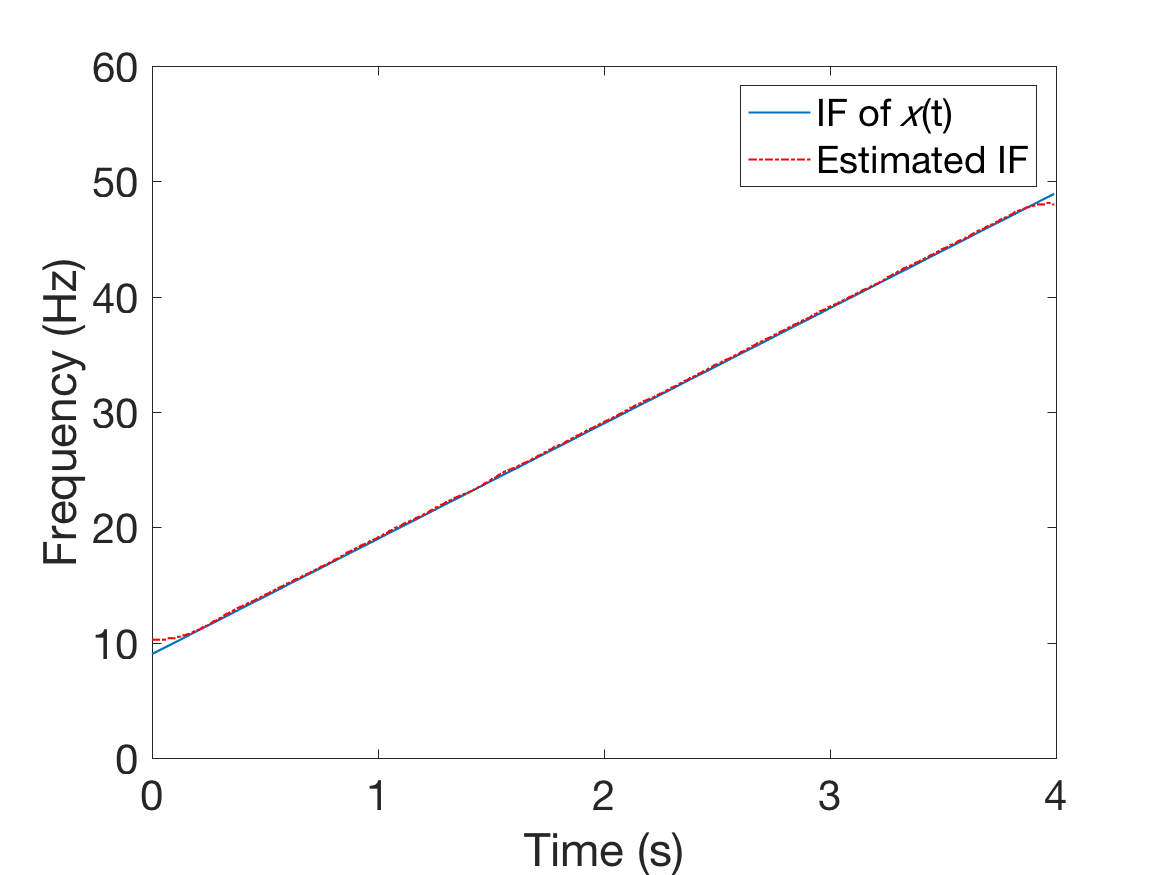}}
\\
\resizebox{3.0in}{2.0in} {\includegraphics{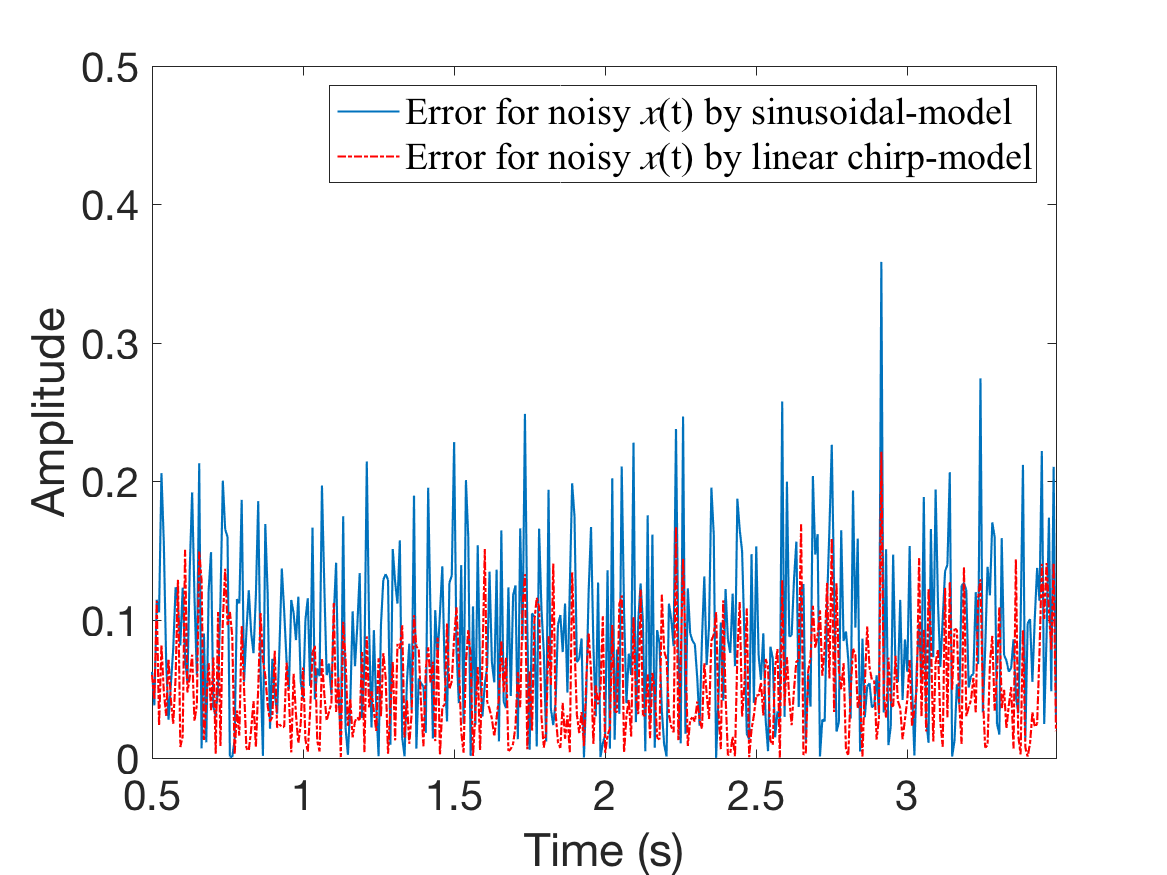}}
	\quad &
	\end{tabular}
	\caption{\small Example of signal $x(t)$ in \eqref{one_chirp} contaminated by a noise with 10dB. 
		Top-left:  Noisy $x(t)$;  
		Top-right:  IF $\phi\rq{}(t) $ of $x(t)$ (solid line) and estimated IF $\wh \eta(t)$ of noisy $x(t)$ (red dot-dashed line); 	Bottom: Absolute recovery errors for noisy $x(t)$ by  sinusoidal signal-based model (blue line) and by linear chirp-based model with an estimated $\phi^{\gp\gp}(t)$ (red dot-dashed line). }
	\label{fig:one_chirp_noise}
\end{figure}
Let $x(t)$ be a linear chirp given by 
\begin{equation}
\label{one_chirp}
x(t)=\cos\big(2\pi(9t+5 t^2)\big), \quad t\in [0, 4),   
\end{equation}
where the number of sampling points is $N=512$ and the sampling rate is 128Hz. 
The IF of  $x(t)$ is  $\phi'(t)=9+10t$, and the chirp rate of $x(t)$ is  $\phi''(t)=10$. 
In the top row of Fig.\ref{fig:one_chirp}, we show  the waveform of $x(t)$ and IF $\phi'(t)$.  In this example and next example, we simply let $\gs(t)\equiv \gs=\frac 1{16}$. The estimated IF, $\big(\wh \eta(t_m)\big)_{ 0\le m<N}$  (where $\wh \eta(t_m)=\wc \eta(t_m)$), is shown on the top-right panel in Fig.\ref{fig:one_chirp} as a red dot-dashed line. The absolute  recovery errors  
$$ 
\big |x(t_m)-2{\rm Re }\big(\wt V_{x}(t_m, \wh \eta(t_m))\big)\big| \quad {\rm and}
\quad \big|x(t_m)-2{\rm Re }\big(\sqrt{1-i2\pi \wt r(t_m)\gs^2} \; \wt V_{x}(t_m, \wc \eta(t_m))\big)\big|
$$ 
by  sinusoidal signal-based model  and by linear chirp-based model with $\wt r(t_m)$ as an approximation to $\phi^{\gp\gp}(t)$ by five-point formula are shown on the bottom-left panel in Fig.\ref{fig:one_chirp} as a blue line and a red line respectively.  Clearly, the recovery error by linear chirp-based model  is much smaller. In the  bottom-right panel in Fig.\ref{fig:one_chirp}, we provide recovery error  $\big|x(t_m)-2{\rm Re }\big(\sqrt{1-i2\pi \phi^{\gp\gp}(t_m)\gs^2} \; \wt V_{x}(t_m, \wc \eta(t_m))\big)\big|$. We know in this case, the error is very small.  

We also consider the performance of models in noisy environment. The results are provided in Fig.\ref{fig:one_chirp_noise}. Again,   linear chirp-based model performs better than  sinusoidal signal-based model.

\begin{figure}[th]
	\centering
	\begin{tabular}{cc}
		\resizebox{3.0in}{2.0in}{\includegraphics{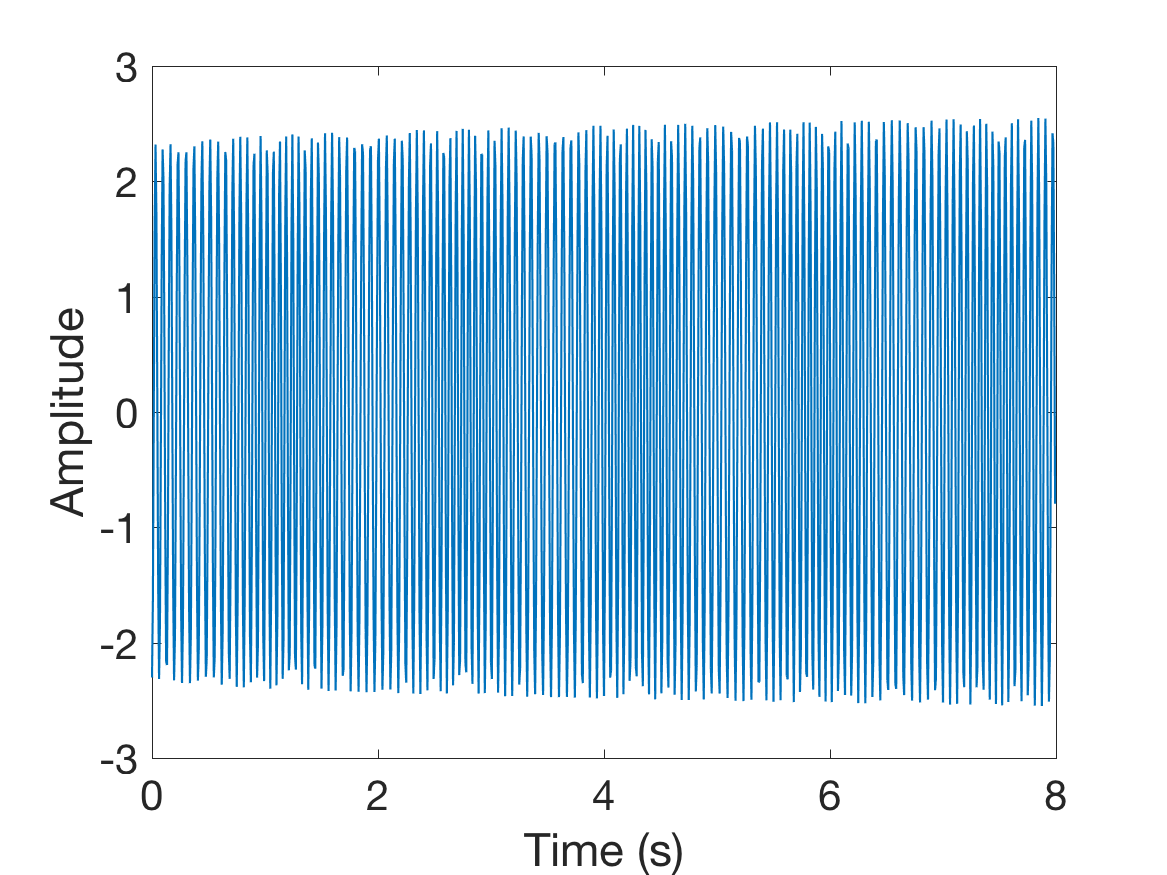}}
		\quad & 
\resizebox{3.0in}{2.0in}{\includegraphics{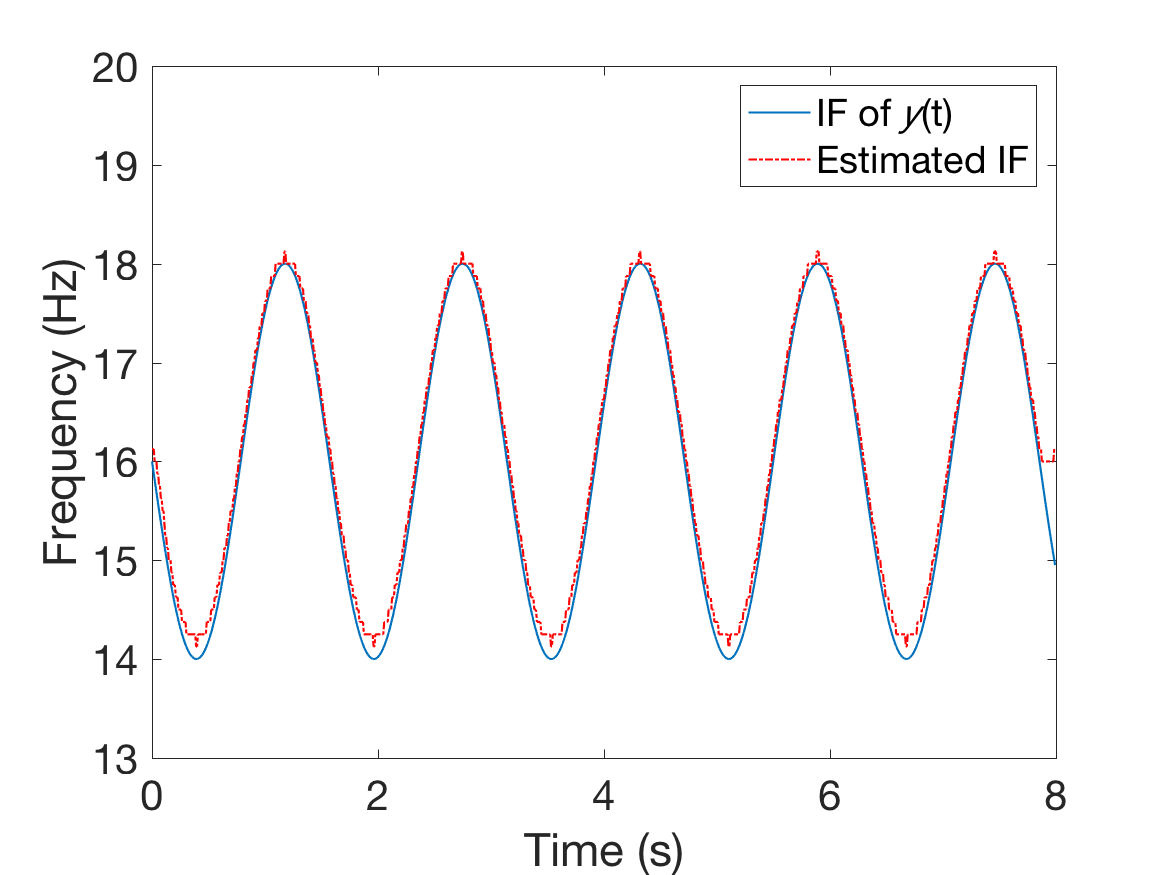}}
\\
 \resizebox{3.0in}{2.0in} {\includegraphics{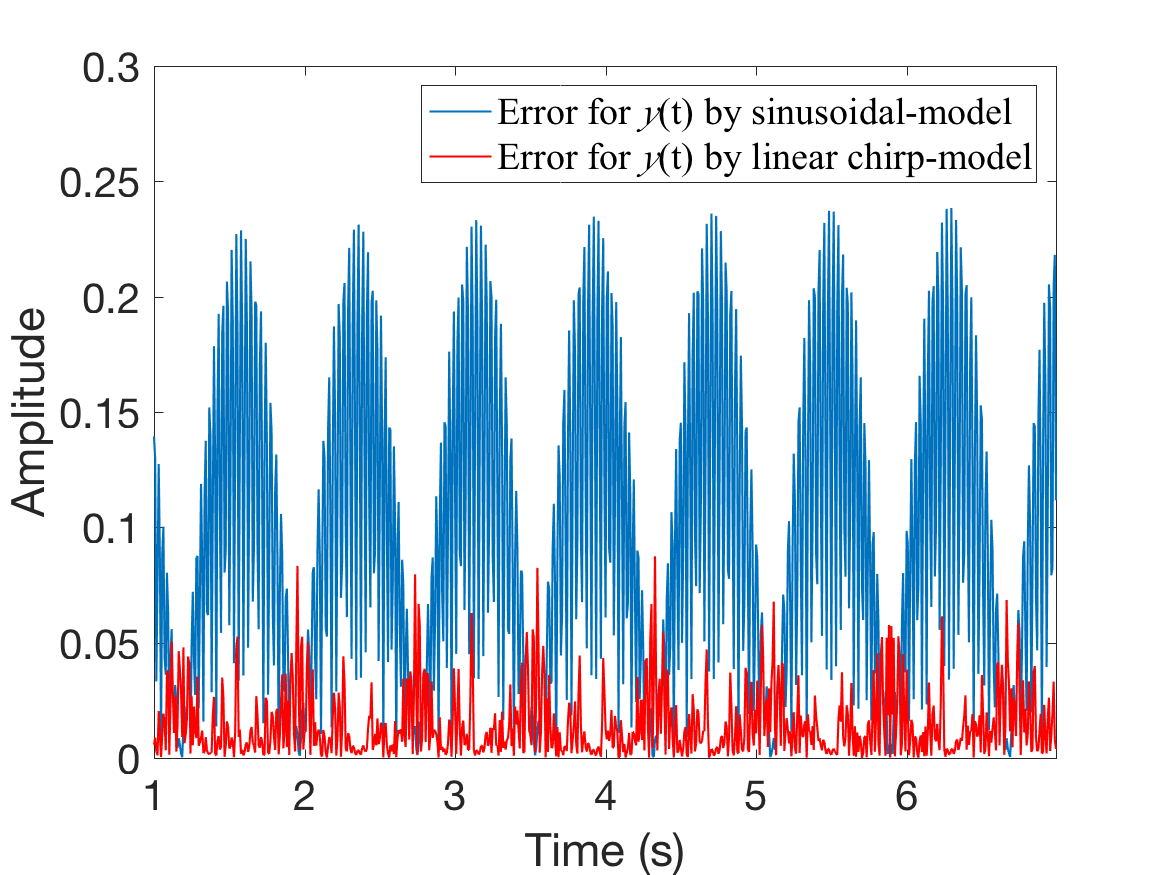}} 
 \quad 
&\resizebox{3.0in}{2.0in} {\includegraphics{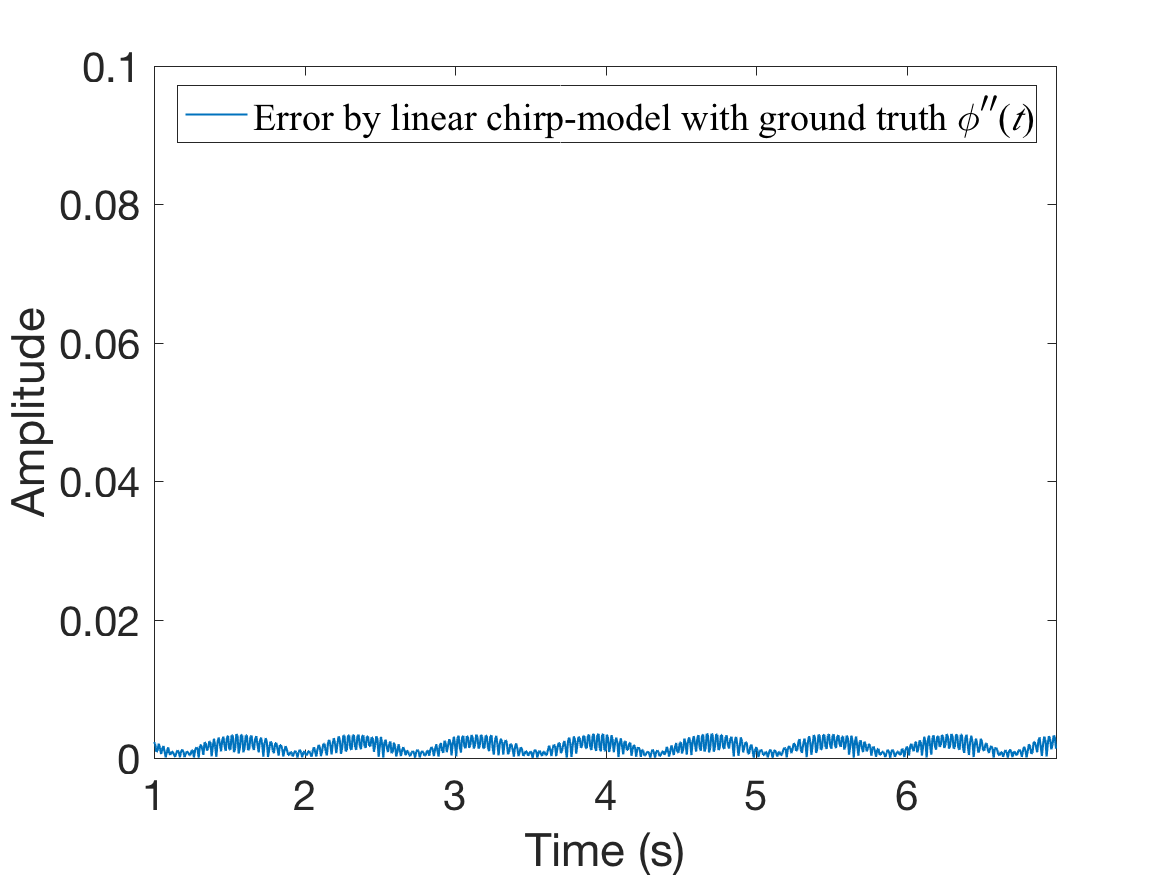}}
	\end{tabular}
	\caption{\small Example of monocomponent signal $y(t)$ in \eqref{one_cosine}. 
		Top-left:  Waveform $y(t)$;  
		Top-right:  IF $\phi\rq{}(t) $ of $y(t)$ (solid line) and estimated IF $\wh \eta(t)$ (red dot-dashed line); 	Bottom-left: Absolute recovery errors for $y(t)$ by  sinusoidal signal-based model (blue line) and by linear chirp-based model with estimated $\phi^{\gp\gp}(t)$ (red line); 
	Bottom-right: Absolute recovery error for $y(t)$ by linear chirp-based model with ground truth $\phi^{\gp\gp}(t)$.}
	\label{fig:one_cosine}
\end{figure}

Next we consider another type of signal given by 
\begin{equation}
\label{one_cosine}
y(t)=\ln(10+\sqrt t)\cos\big(2\pi(16t+0.5\cos(4t))\big), \quad t\in [0, 8).  
\end{equation}
This time IF of $y(t)$ is a  sinusoidal function. $y(t)$ is uniformly sampled with the number of sampling points $N=1024$ and the sampling rate 128Hz.   
The waveform of $y(t)$ and its IF $\phi'(t)$ are shown in the top row of Fig.\ref{fig:one_cosine}, where  the estimated IF $\big(\wh \eta(t_m)\big)_{ 0\le m<N}$ is also provided (on the top-right panel) as a red dot-dashed line. The absolute  recovery errors  by  sinusoidal signal-based model  and by linear chirp-based model with $\wt r(t_m)$ as an approximation to $\phi^{\gp\gp}(t)$ by five-point formula are shown on the bottom-left panel in Fig.\ref{fig:one_cosine} as a blue line and a red line respectively.  By the way, we also provide the recovery error by linear chirp-based model with ground truth $\phi^{\gp\gp}(t)$. In addition, we also provide recovery errors 
for noisy $y(t)$ in Fig.\ref{fig:one_cosine_noise}. Again,  linear chirp-based model results in smaller recovery errors. 

\begin{figure}[th]
	\centering
	\begin{tabular}{cc}
		\resizebox{3.0in}{2.0in}
{\includegraphics{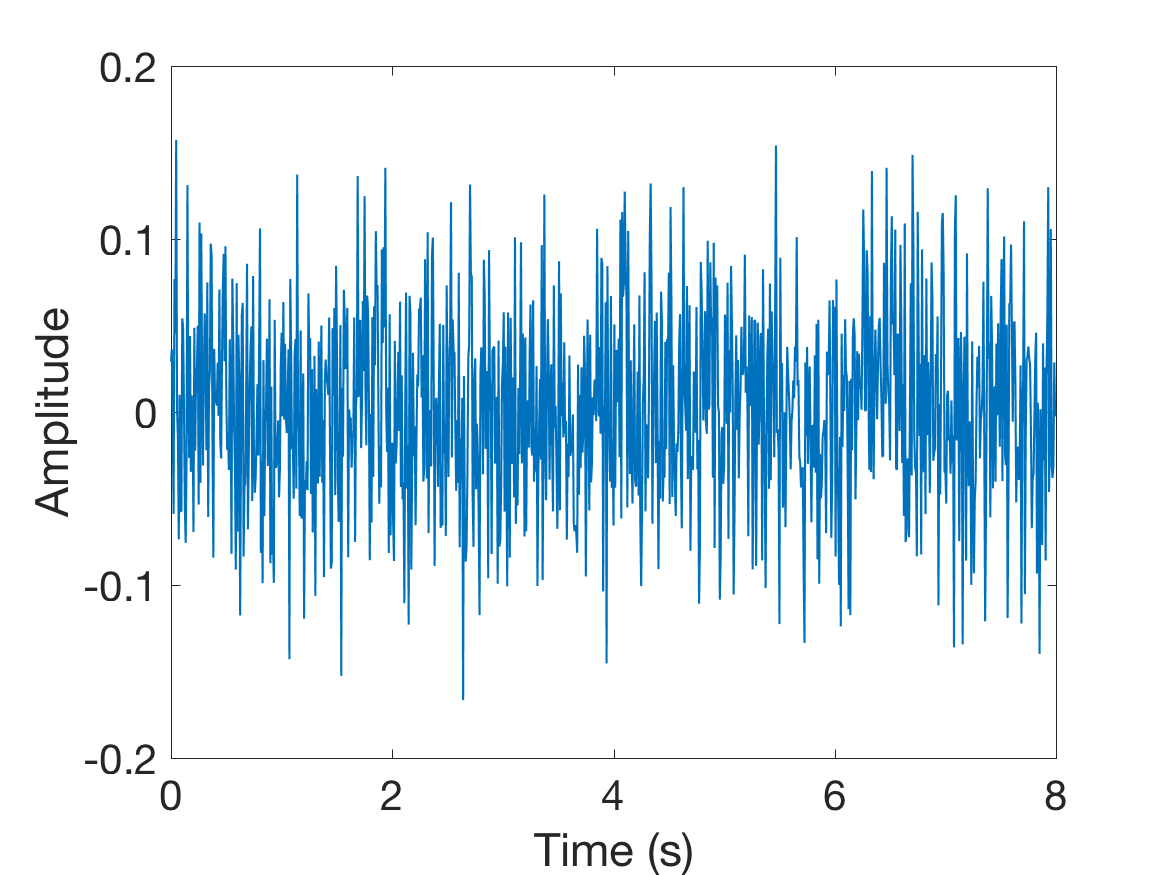}}\quad 
	\quad &\resizebox{3.0in}{2.0in}{\includegraphics{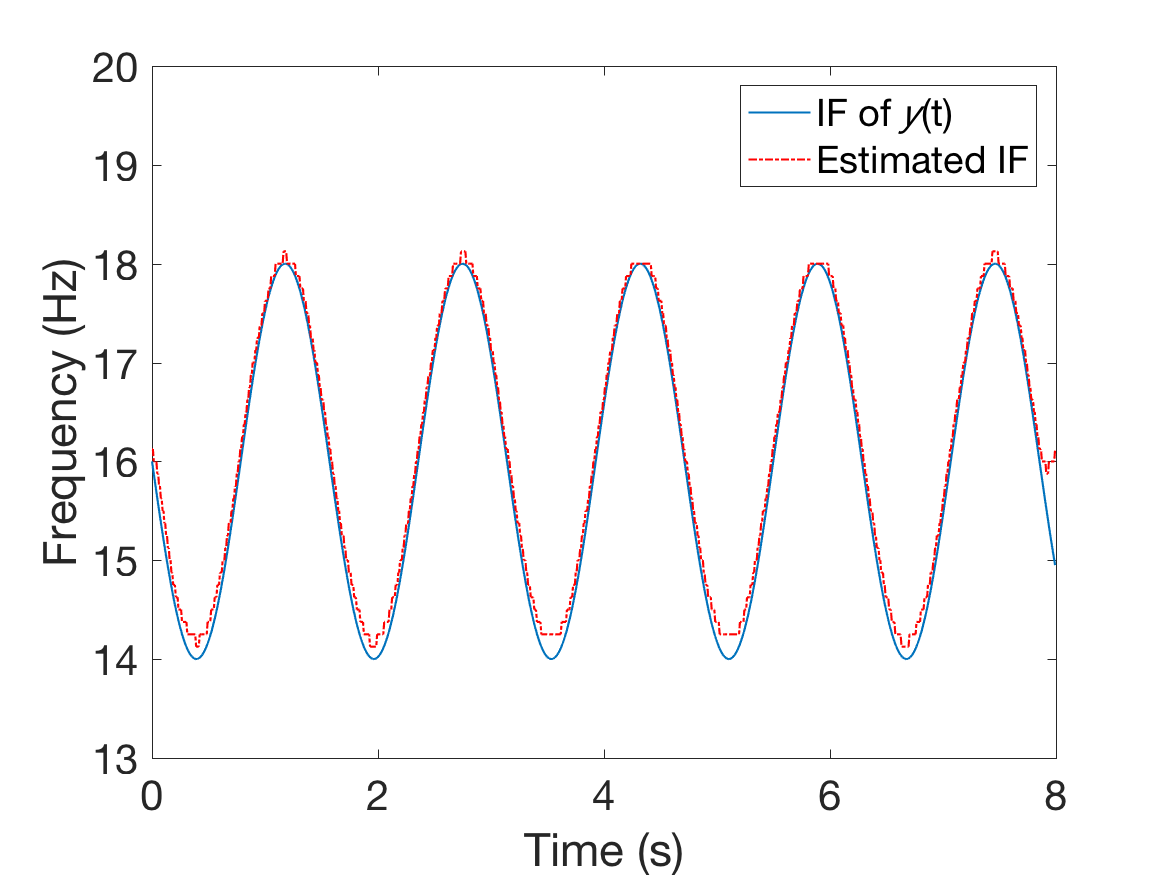}}
\\
\resizebox{3.0in}{2.0in} {\includegraphics{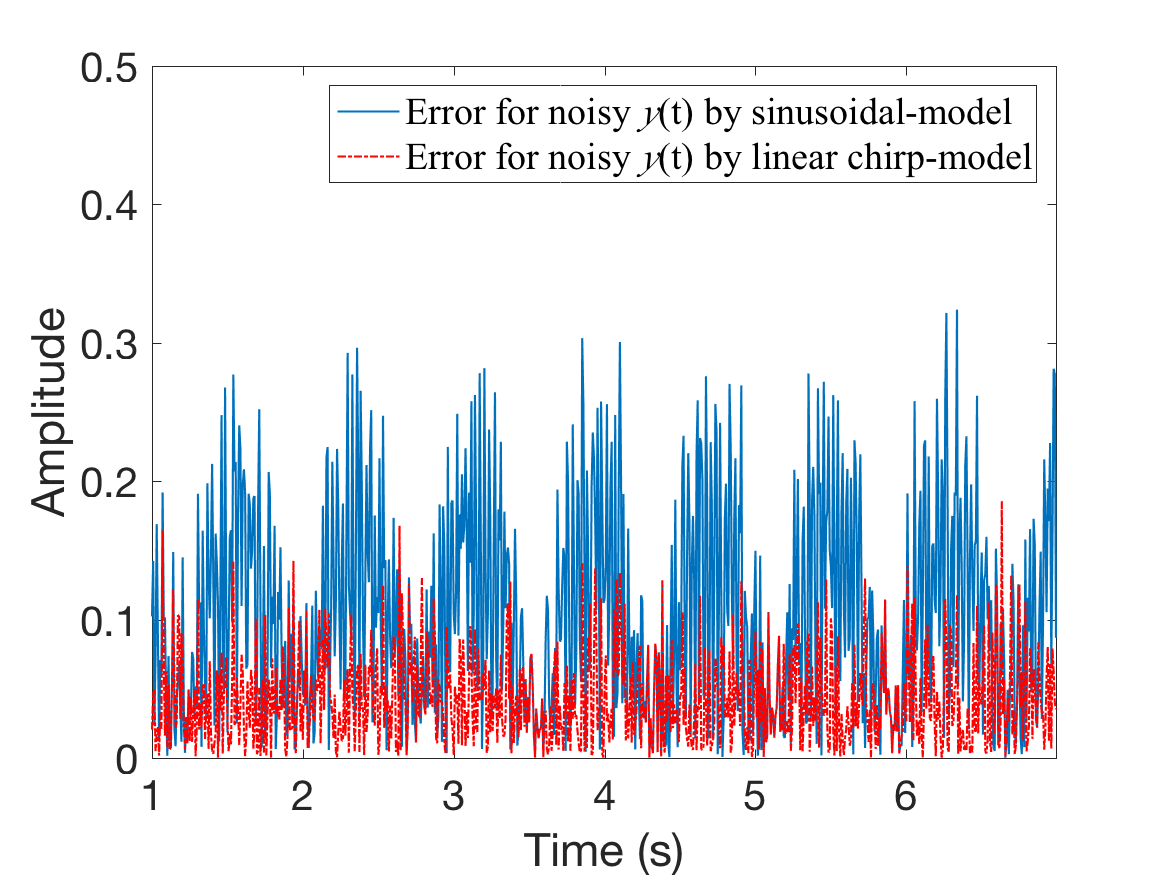}}
	\quad &
	\end{tabular}
	\caption{\small Example of signal $y(t)$ in \eqref{one_cosine} contaminated by a noise with 15dB. 
		Top-left:  Added noise;  
		Top-right:  IF $\phi\rq{}(t) $ of $y(t)$ (solid line) and estimated IF $\wh \eta(t)$ of noisy $y(t)$ (red dot-dashed line); 	Bottom: Absolute recovery errors for noisy $y(t)$ by  sinusoidal signal-based model (blue line) and by linear chirp-based model with an estimated $\phi^{\gp\gp}(t)$ (red dot-dashed line). }
	\label{fig:one_cosine_noise}
\end{figure}

\bigskip 

\begin{figure}[th]
	\centering
	\begin{tabular}{cc}
		\resizebox{3.0in}{2.0in}{\includegraphics{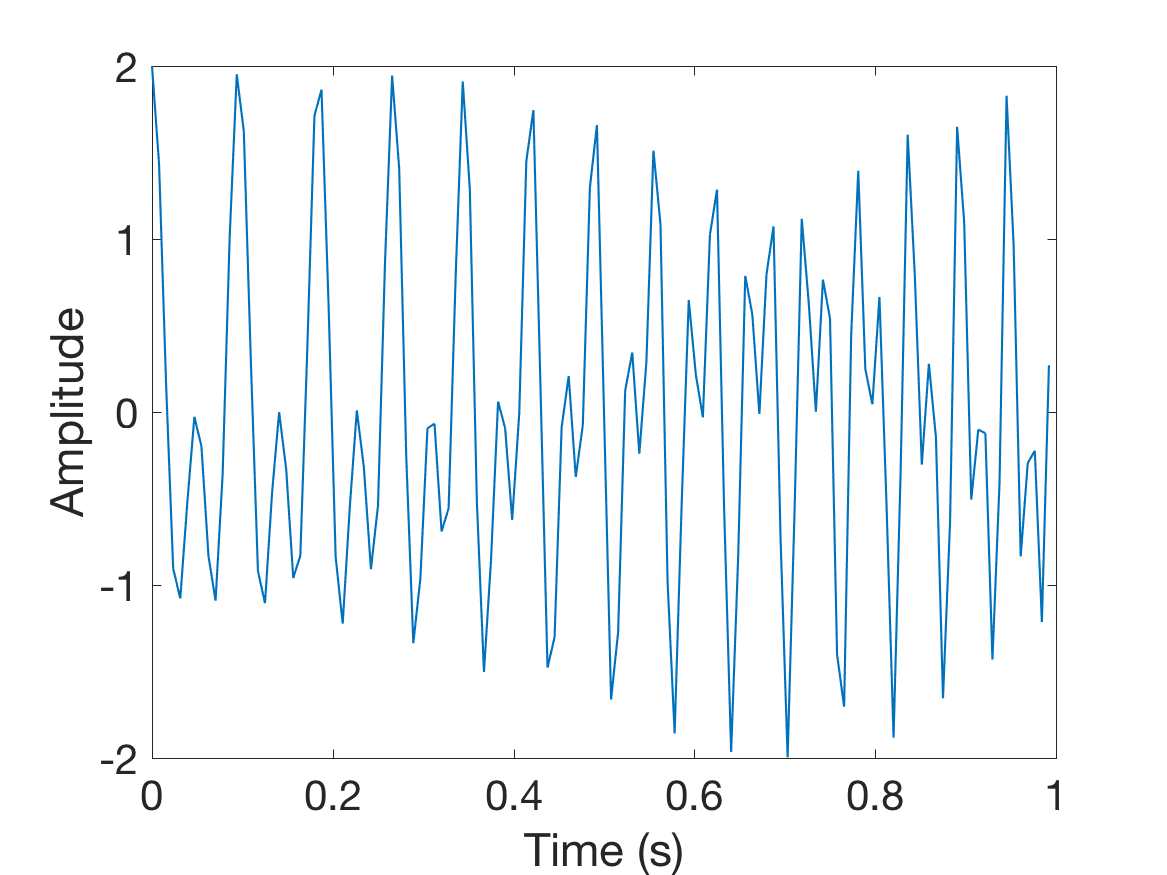}}
		\quad & 
\resizebox{3.0in}{2.0in}{\includegraphics{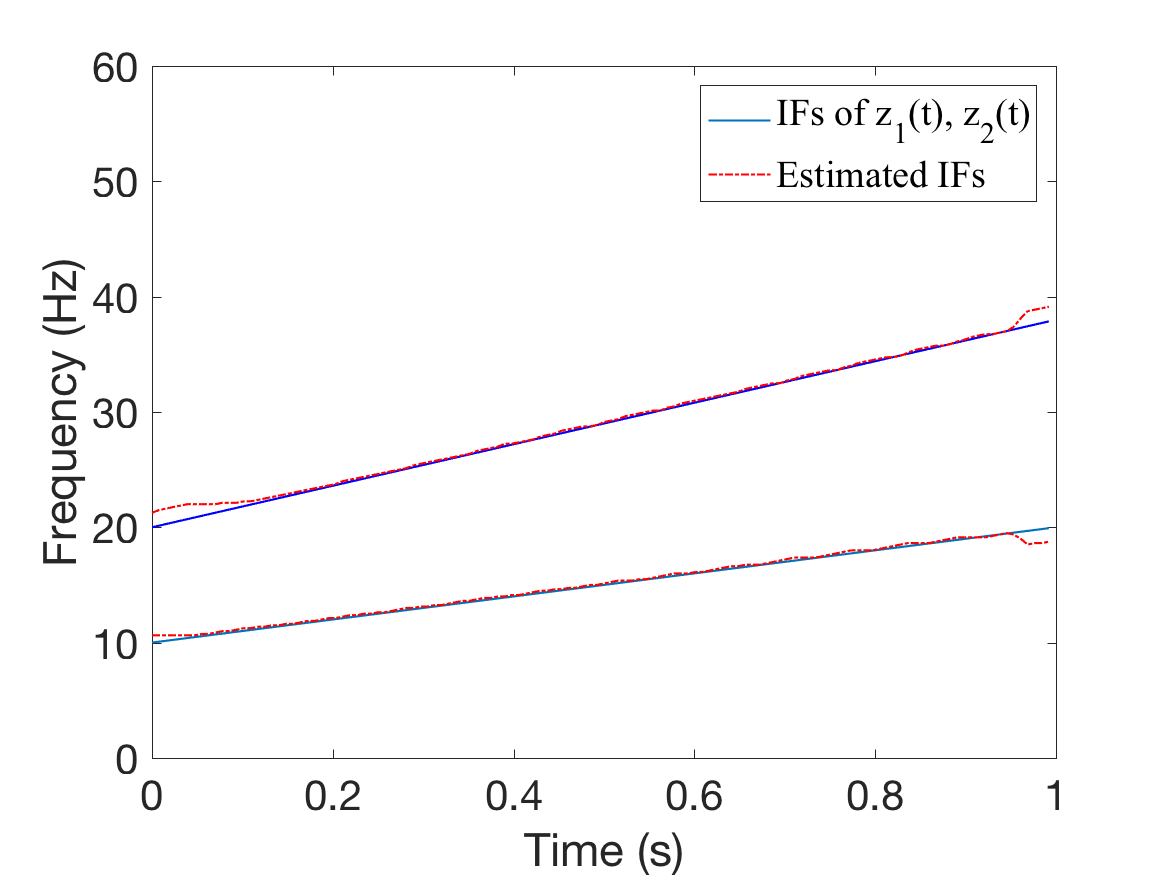}}
\\
\resizebox{3.0in}{2.0in}
{\includegraphics{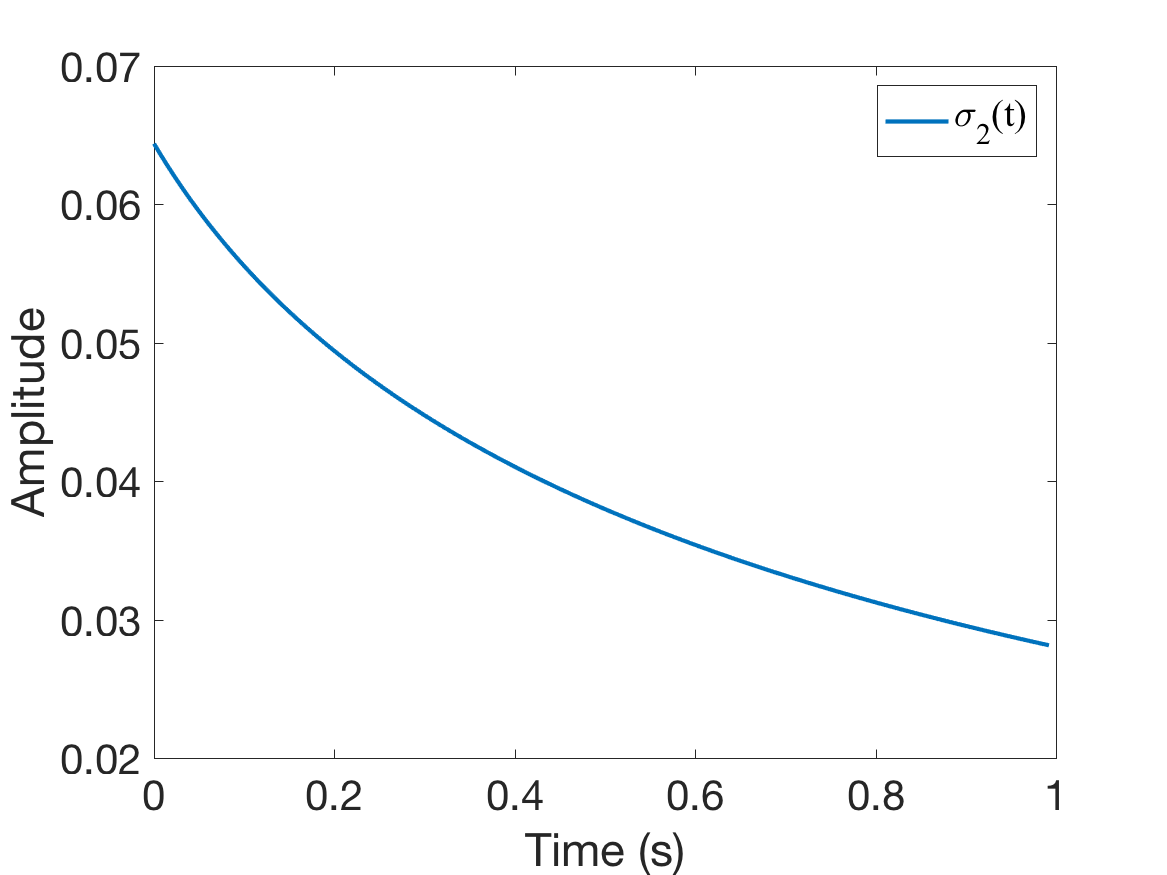}}
		\quad & \\
 \resizebox{3.0in}{2.0in} {\includegraphics{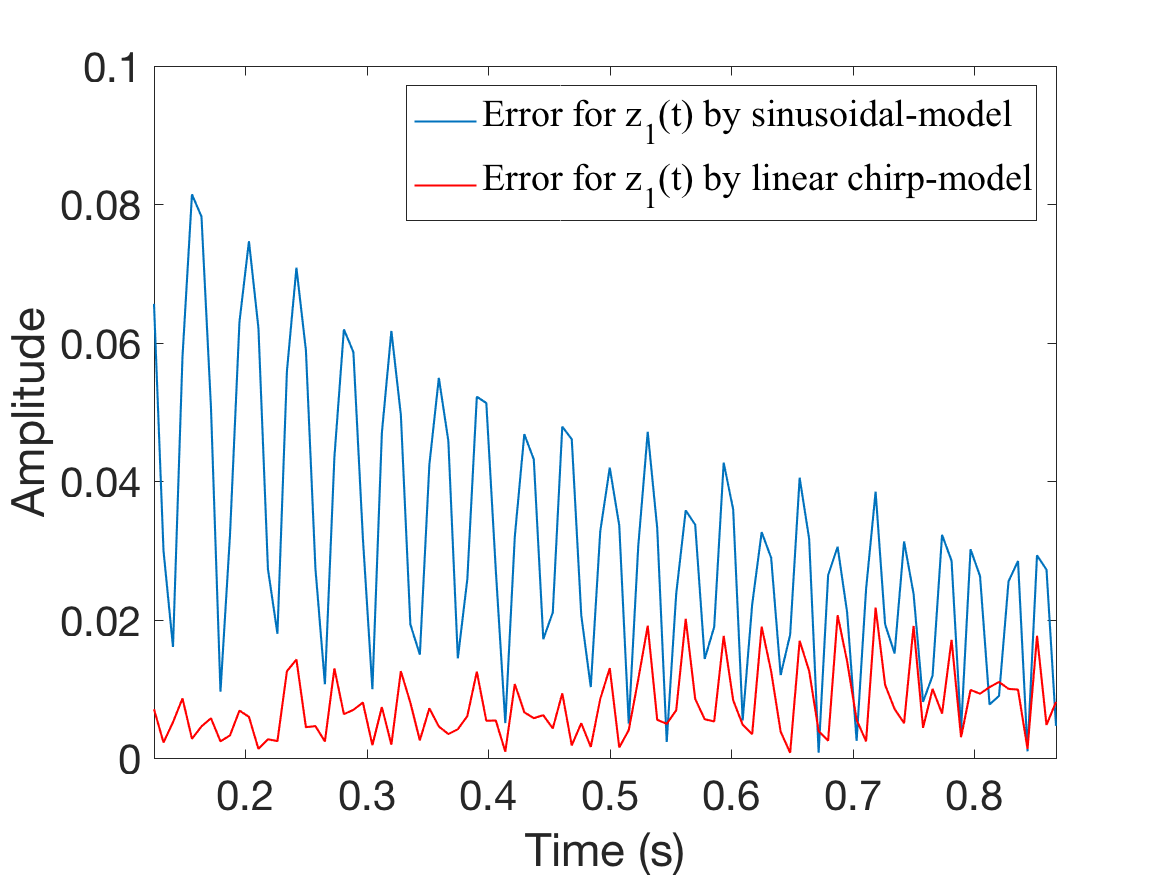}} 
 \quad 
&\resizebox{3.0in}{2.0in} {\includegraphics{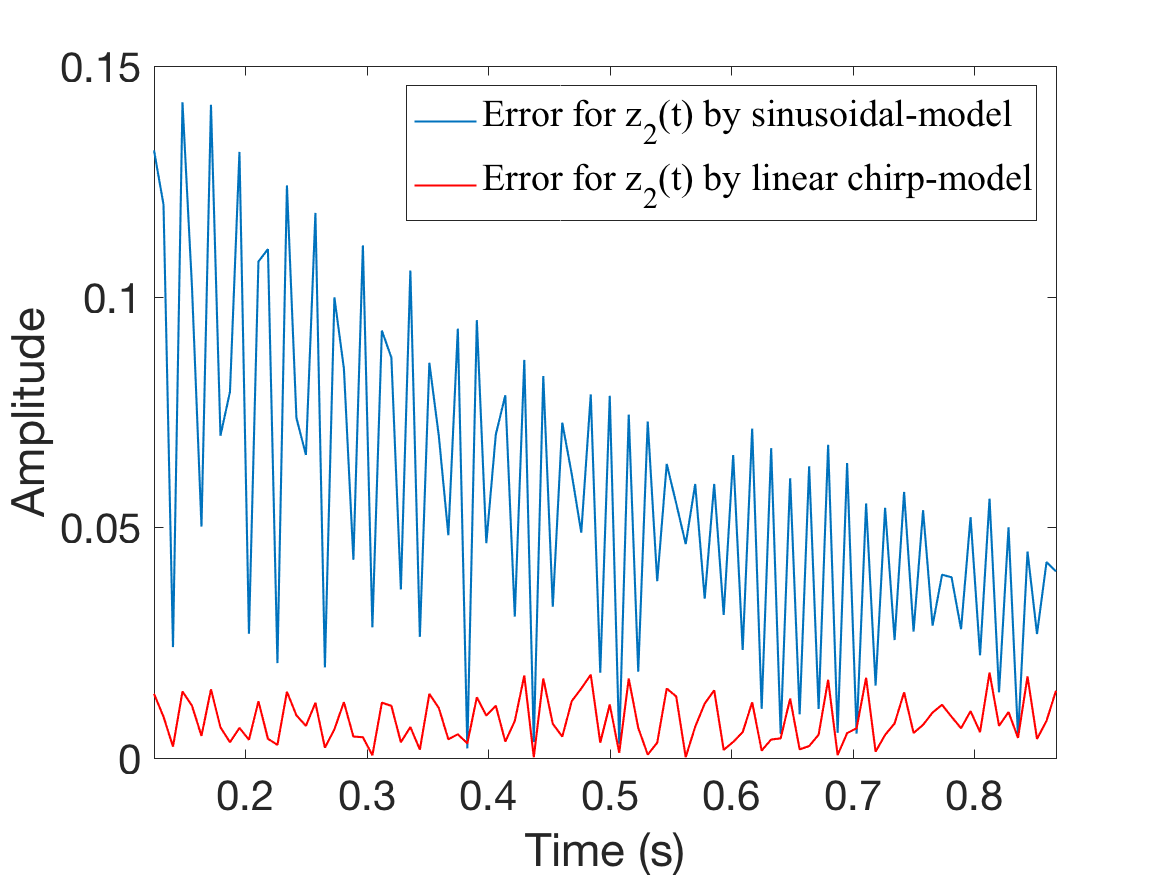}}
	\end{tabular}
	\caption{\small Example of monocomponent signal $z(t)$ in \eqref{two_LFM}. 
		Top-left:  Waveform $Z(t)$;  
		Top-right:  IFs $\phi_1\rq{}(t), \phi_2\rq{}(t)$ of $z_1(t), z_2(t)$ (solid line) and estimated IFs $\wh \eta_1(t), \wh \eta_2(t)$ (red dot-dashed line); 
		Middle: Time-varying parameter $\gs_2(t)$;
				Bottom-left: Absolute recovery errors for $z_1(t)$ by  sinusoidal signal-based model (blue line) and by linear chirp-based model with estimated $\phi_1^{\gp\gp}(t)$ (red line); 
	Absolute recovery errors for $z_2(t)$ by  sinusoidal signal-based model (blue line) and by linear chirp-based model with estimated $\phi_2^{\gp\gp}(t)$ (red line).}
	\label{fig:two_chirps}
\end{figure}

Finally, we consider a two-component signal given by 
\begin{equation}
\label{two_LFM}
z(t)=z_1(t)+z_2(t)=\cos\big(2\pi(10t+5t^2)\big)+ \cos\big(2\pi(20t+9t^2)\big), \quad t\in [0, 1),   
\end{equation}
where $z(t)$ is uniformly sampled with the number of sampling points $N=128$ and the sampling rate 128Hz.   
The waveform of $z(t)$ and its IFs $\phi'_1(t), \phi'_2(t)$ are shown in the top row of Fig.\ref{fig:two_chirps}. 
In this example, we also use a time-varying parameter $\gs(t)$. The authors in \cite{LCHJJ18} propose an algorithm to select $\gs(t)$. Here we let $\gs(t)=\gs_2(t)$ suggested in \cite{LCHJJ18}. $\gs_2(t)$ is shown in the second row of  Fig.\ref{fig:two_chirps}. 
 With this $\gs_2(t)$, the estimated IFs $\big(\wh \eta_\ell(t_m)\big)_{ 0\le m<N}, \ell=1, 2$ are provided on the top-right panel of Fig.\ref{fig:two_chirps} as red dot-dashed lines. The absolute recovery errors for $z_1(t)$ by  sinusoidal signal-based model  and by linear chirp-based model with $\wt r_1(t_m)$ as an approximation to $\phi^{\gp\gp}_1(t)$ by five-point formula are shown on the bottom-left panel of Fig.\ref{fig:two_chirps} as a blue line and a red line respectively, while the recovery errors for $z_2(t)$ are presented on the bottom-right panel of Fig.\ref{fig:two_chirps}. Linear chirp-based model clearly separates the components $z_1(t), z_2(t)$ of $z(t)$ better than  sinusoidal signal-based model. 
 
 \clearpage 
 
 We also consider the recovery results with $\gs(t)\equiv \gs=\frac 1{16}$. Here instead of providing pictures, we use the relative root of mean square error (RMSE) to evaluate the error of component recovery.
 The RMSE is defined by 
\begin{equation*}
{\rm RMSE}_\upsilon := \frac{1}{K} \sum_{k=1}^{K} \frac{\|\upsilon_k-\wh \upsilon_k\|_2} {\|\upsilon_k\|_2}, 
\end{equation*} 
where $\upsilon$ is a vector and $\wh\upsilon$ is an estimation of  $\upsilon$. Note that as mentioned above,  
the RMSE is calculated over $[N/8+1:7N/8]$ 
to avoid the errors near the end points, which are usually large.  We use ASTFT${}_{\rm si}$ and ASTFT${}_{\rm lc}$ to denote the  sinusoidal signal-based model and the linear chirp-based model respectively. 
 From Table \ref{tab1}, we see the recovery formula based on adaptive STFT with time-varying $\gs(t)$ results in smaller component recovery errors than that based on the conventional STFT. 
  
\begin{table}
\begin{center}
\begin{tabular}{|l|c|c|c|c|c}
\hline
&{\small ASTFT${}_{\rm si}$} 
&{\small ASTFT${}_{\rm lc}$ }
&{\small ASTFT${}_{\rm lc}$ with $\phi^{\gp\gp}(t)$ }
 \\\hline
$\gs=\frac1{16}$ &0.3254&0.0308&0.0123\\\hline
$\gs=\gs_2(t)$&0.1458 &0.0264&0.0090\\\hline
\end{tabular}
\caption{{\small RMSEs for $z(t)$.} }
\label{tab1}
\end{center}
\end{table}



\end{document}